\documentclass[10pt, a4paper]{article}
\usepackage{a4wide}
\usepackage[pdftex]{graphicx}
\usepackage[dvipdfmx]{color}
\usepackage{color}
\usepackage{bm}
\usepackage{amsmath}
\usepackage{amssymb}
\usepackage{amsfonts}
\usepackage{amsmath}
\usepackage{ascmac}
\usepackage[english]{babel}
\usepackage{mathrsfs}
\title{A construction of Exceptional Weyl Group $W(F_4)$ and $W(E_8)$ using Quaternion, and the lattice in 16-dimensional Euclidean space}
\author{Misaki Ohta,University of the Ryukyus \thanks{Department of physics, Email: e193225@eve.u-ryukyu.ac.jp or apple.designed@icloud.com}}
\begin{document}
\maketitle

\begin{abstract}
It is mentioned that there is a subalgebra isomorphic to the alternating group $2 \cdot A_4$ as a subalgebra of the Quaternion over integers and half-integers called Hurwitz quaternionic integers $\mathscr{H}$ in the book by J.H.Conway and Neil J. A. Sloane. In this paper, I have followed this book and extended Quaternion over integers and half-integers to have duality, and proved that a subalgebra in it isomorphic to Exceptional Weyl group $W(F_4)$. I have also found a method of constructing the $16$-dimensional lattice $\Lambda_{16}$ which seems to be isomorphic to the lattice called the Barnes-Wall lattice $\Lambda_{\text {Barnes-Wall }}$, which is currently considered to be very dense (although this remains to be discussed) using the Dual Quaternion. Lastly, I briefly mention how to construct an exceptional Weyl group $W(E_8)$ using an Octonion and Dual Quaternion.
\end{abstract}

\section{Introduction}
\ In the history of physics, beautiful symmetries have been found arising from the algebraic structures they satisfy, such as the Pauli matrices of Quantum Mechanics and the Gamma matrices of Dirac in relativistic quantum mechanics. And there is always some inevitability behind the appearance of symmetries. It is known that in order to understand the structure of symmetry, it is necessary to consider the structure of the symmetry (i.e., the group) itself and the object on which it acts in depath, and such studies have actually been conducted. One of the symmetries that will play a leading role in this paper is called Quaternions, which is useful in representing the rotation of coordinates, discovered by Hamilton on the $16$th of October $1843$. It is frequently used not only in a physical context but also in a mathematical context and has many applications. The fact that it is widely applied means that it plays a very profound core role in this world. On the other hand, in terms of the objects on which the group acts, there is a called Root System, which will be the main subject of this paper, closely related to Lie algebras and Lie groups, and essential from the point of view of describing the most advanced theoretical physics of our time. Root System is simply a set of vectors generated by reflection, but it is generally known that it can be defined even in high order coordinates, and over a long period of history, the classification of the types of reflection was completed. Among the classified systems, those that deviate from the general classification are called exceptional Root Systems, such as $G_2,F_4$ and $E_8$. In particular, the $E_8$ root system in $8$th dimension is considered to be the largest among the root systems of exceptional type, and this also appears as a necessity in the world of theoretical physics, the $E_8$ lattice forms an 8-dimensional dense packing as well. In addition, The number of vectors that are the shortest of a lattice, i.e., the number of the densest packing, tells us how reasonable and natural the lattice is, because the greatest number of shortest vectors means that they are most evenly distributed on the sphere of the space. And for vectors that form a dense packing, it is always possible to discuss orthogonal transformations that transfer between those vectors, i.e., it is a symmetry discussion for that lattice. The orthogonal groups acting on Root System are called Weyl groups. For example, the Weyl groups of $G_2,F_4$ and $E_8$ are written as $W(G_2),W(F_4)$ and $W(E_8)$, respectively. In this paper, I will first give a duality to Quaternion, and then describe the existence of a subalgebra in the algebraic system corresponding to the reflection of exceptional type $W(F_{4})$ in 4-dimensional spaces described above. Next, I will briefly study a lattice, which can be constructed by this subalgebra, with $4320$ shortest vectors similar (or isomorphic) to the Barnes-Wall lattice, which is said to be the close-packed of $16$-dimensional Euclidean space through the homomorphism to the code(vector space) over $\mathbb{F}_{2}$ which is invariant by the action of the group $2^{4} \cdot A_{8}$, which is the maximal subgroup of Mathieu group $M_{24}$ called 'Miracle Octad Group'. Lastly, I briefly mention how to construct an exceptional Weyl group $W(E_8)$ using an Octonions and Dual Quaternion.

\subsection{Quaternion and Hurwitz quaternionic integers}
Let $\mathbb{C}$ donate the field of complex numbers $\{z = x+iy|x,y \in \mathbb{R}\}$, where $i^2 = -1$, and let $\mathbb{H}$ donate the filed of Quaternions $\{z = x_0 + x_1i + x_2j + x_3k |x_i \in \mathbb{R}\}$ where
$$
\begin{array}{l}
i^{2}=j^{2}=k^{2}=i j k=-1 \\
i j=-j i=k, \quad j k=-k j=i, \quad k i=-i k=j
\end{array}
$$
The conjugate $\bar{z}$ of $z$ is defined to be
$$
\begin{array}{ll}
\bar{z}=x-i y & \in \mathbb{C} \\
\bar{z}=x_{0}-i x_{1}-j x_{2}-k x_{3} & \in \mathbb{H}
\end{array}
$$
and then the norm $N:z \in \mathbb{C},\mathbb{H} \mapsto z\bar{z} \in \mathbb{R}$ of $z$ is
$$
\begin{array}{l}
x^{2}+y^{2} \\
x_{0}^{2}+x_{1}^{2}+x_{2}^{2}+x_{3}^{2}
\end{array}
$$
respectively.There are the Gaussian and Eisenstein integers
$$
\begin{array}{l}
\mathscr{G}=\{a+i b: a, b \in \mathbb{Z}\} \subset \mathbb{C}, \\
\mathscr{E}=\{a+\omega b: a, b \in \mathbb{Z}\} \subset \mathbb{C},
\end{array}
$$

\begin{figure}[h]
  \begin{center}
 \includegraphics[keepaspectratio, scale=0.1]
      {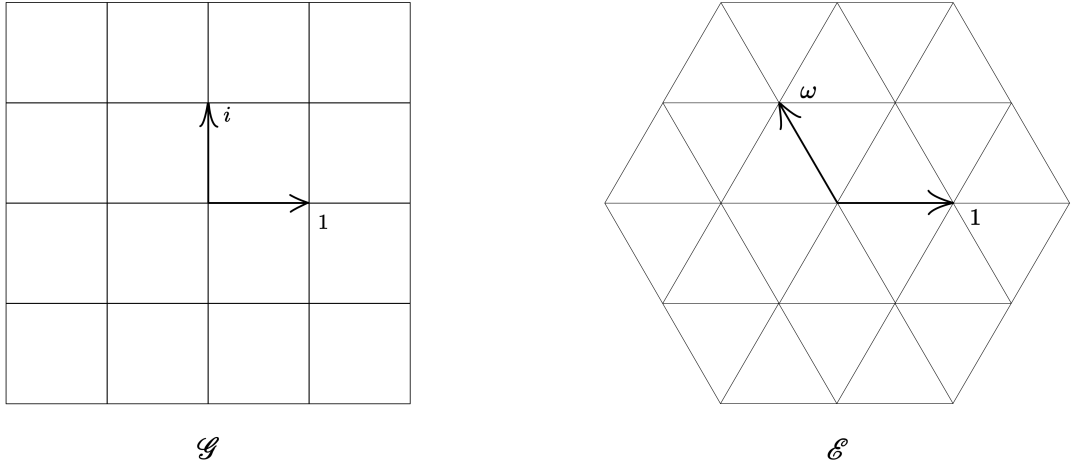}
 \caption{Gaussian and Eisenstein integers}
\label{}
\end{center}
\end{figure}

where $\omega = (-1 + i\sqrt{3})/2$, and the Hurwitz Quaternionic integers
$$
\mathscr{H}=\{a+i b+j c+k d: a, b, c, d \in \mathbb{Z} \text { or } a, b, c, d \in \mathbb{Z}+1 / 2\} \subset \mathbb{H}.
$$
According to "Sphere Packings, Lattices and Groups" by John Conway, Neil J. A. Sloane[1,P53], the set $\mathcal{S} \subset \mathscr{H} \subset \mathbb{H}$ generated by the multiplication of the following generators forms a group
$$
i, j, k, \frac{1}{2}(-1+i+j+k)  \in \mathbb{H}
$$
and there is isomorphism
$$
\mathcal{S} \cong 2 \cdot A_{4}
$$
where $A_4$ is the alternating group of degree $4$[See Appendix 5.4]. The fields $\mathscr{G}, \mathscr{E}$ and $\mathscr{H}$ can also be interpreted as the group rings over $\mathbb{Z}$ or $\mathbb{Z} + 1/2$ of the following groups  (may divided by the center)
$$
\begin{array}{l}
G_{1}=< i>\cong \mathbb{Z} / 4 \mathbb{Z} \\
G_{2}=<\omega> \cong \mathbb{Z} / 3 \mathbb{Z} \\
G_{3}=\langle i, j\rangle \cong Q_{8}=\left\langle a, b \mid a^{4}=e, a^{2}=b^{2}, b a=a^{-1} b\right\rangle
\end{array}
$$
In general, an element $x$ of group ring $R[G]$ of $G$ can be written as
$$
x=\sum_{g_{i} \in G} g_{i} x_{i} \quad (x_i \in \mathbb{R})
$$
if I follow the definition of conjugation above, then the conjugation $\bar{x}$ of $x$ should be
$$
\bar{x}=\sum_{g_{i} \in G} g^{-1}_{i} x_{i} \quad\left(x_{i} \in \mathbb{R}\right)
$$
and if the group is an orthogonal group($G \subset O(n)$), then it is simply shown that $x$ is Hermitian $\bar{x}={}^{t} x(x=x^{\dagger})$. Because
$$
\begin{aligned}
\bar{x} &=\sum_{g_{i} \in G} x_{i} g_{i}^{-1}=\sum_{g_{i} \in G} x_{i}^{t} g_{i} \\
&={ }^{t}\left(\sum_{g_{i} \in G} x_{i} g_{i}\right)={ }^{t} x
\end{aligned}
$$
and when $N(x) = x\bar{x} = 1$, in other words ${ }^{t} x=x^{-1}$, it means
$$
x \in O(n) \quad (x \in R[G] \quad G \subset O(n))
$$
And on the contrary, if $G$ is an orthogonal group and $x \in R[G]$ is also an enelent of an orthogonal group, then the norm defined here is determined to be $1$.
\subsection{The representation of quaternions}
Just as the algebraic structure of Pauli matrices in quantum mechanics
$$
\begin{aligned}
&\sigma_{1}^{2}=\sigma_{2}^{2}=\sigma_{3}^{2}=1 \\
&\sigma_{1} \sigma_{2}=-\sigma_{2} \sigma_{1}=i \sigma_{3} \\
&\sigma_{2} \sigma_{3}=-\sigma_{3} \sigma_{2}=i \sigma_{1} \\
&\sigma_{3} \sigma_{1}=-\sigma_{1} \sigma_{3}=i \sigma_{2}
\end{aligned}
$$
allows the matrix to be determined
$$
\sigma_{}=\left(\begin{array}{ll}
1 &  \\
 & 1
\end{array}\right), \quad \sigma_{1}=\left(\begin{array}{ll}
 & 1 \\
1 &
\end{array}\right), \quad \sigma_{2}=\left(\begin{array}{cc}
 & -i \\
i &
\end{array}\right), \quad \sigma_{3}=\left(\begin{array}{cc}
1 &  \\
 & -1
\end{array}\right)
$$
or the algebraic structure of Dirac spinor in QFT
$$
\begin{aligned}
\left(\gamma_{0}\right)^{2} &=1,\left(\gamma_{j}\right)^{2} =-1 \quad\left(j=1,2,3\right)\\
\left\{\gamma_{\mu}, \gamma_{\nu}\right\} &=\gamma_{\mu} \gamma_{\nu}+\gamma_{\nu} \gamma_{\mu}=0 \quad(\mu \neq \nu)
\end{aligned}
$$
allows the matrix to be determined
$$
\begin{array}{l}
\gamma_{0}=\left(\begin{array}{ll}
 & E \\
E &
\end{array}\right), \quad \gamma_{j}=\left(\begin{array}{cc}
 & \tau_{j} \\
-\tau_{j} &
\end{array}\right) \\
\text{where  }\quad\tau_{1}=\left(\begin{array}{ll}
 & 1 \\
1 &
\end{array}\right), \quad \tau_{2}=\left(\begin{array}{cc}
 & -i \\
i &
\end{array}\right), \quad \tau_{3}=\left(\begin{array}{cc}
1 &  \\
 & -1
\end{array}\right)
\end{array}
$$
from the stracture of Quaternions
$$
\begin{aligned}
&i^{2}=j^{2}=k^{2}=i j k=-1 \\
&i j=-j i=k, \quad j k=-k j=i, \quad k i=-i k=j
\end{aligned}
$$
it allows the matrix to be determined using isomorphism $\rho: <i,j> \rightarrow O(4,\mathbb{R})$
$$
\begin{array}{l}
 (i)  \rho=\left(\begin{array}{cccc}
 & 1 &  &  \\
-1 &  &  &  \\
 &  &  & 1 \\
 &  & -1 &
\end{array}\right) \\\\
 (j)  \rho=\left(\begin{array}{cccc}
 &  & 1 &  \\
 &  &  & -1 \\
-1 &  &  &  \\
 & 1 &  &
\end{array}\right) \\\\
 (k)  \rho=\left(\begin{array}{cccc}
 &  &  & -1 \\
 &  & -1 &  \\
 & 1 &  &  \\
1 &  &  &
\end{array}\right)
\end{array}
$$
Using this representation, the group $\mathcal{S} \subset \mathscr{H} \subset \mathbb{H}$ described above yields
$$
\mathcal{S}\cong
\left\langle\left(\begin{array}{rrrr}
 & 1 &  &  \\
-1 &  &  &  \\
 &  &  & 1 \\
 &  & -1 &
\end{array}\right),\left(\begin{array}{rrrr}
 &  & 1 &  \\
 &  &  & -1 \\
-1 &  &  &  \\
 & 1 &  &
\end{array}\right) \frac{1}{2} \left(\begin{array}{rrrr}
-1 & 1 & 1 & -1 \\
-1 & -1 & -1 & -1 \\
-1 & 1 & -1 & 1 \\
1 & 1 & -1 & -1
\end{array}\right)\right\rangle \cong 2 \cdot A_4
$$

\subsection{Root System}
\begin{figure}[h]
  \begin{center}
 \includegraphics[keepaspectratio, scale=0.3]
      {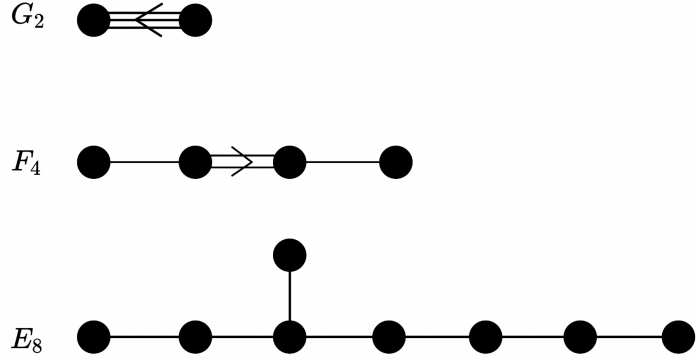}
 \caption{Dynkin diagrams of exceptional root system}
\label{}
\end{center}
\end{figure}

\subsubsection{Root system $F_4$}
In this paper, I use the concept of root system here.\\\\
Let $\left\{\varepsilon_{1}, \varepsilon_{2}, \varepsilon_{3}, \varepsilon_{4}\right\}$ be the standard basis of
$\mathbb{R}^4$. I set the notation for $F_4$ as generated by the basis. $\varepsilon_{2}-\varepsilon_{3}, \varepsilon_{3}-\varepsilon_{4}, \varepsilon_{4},\left(\varepsilon_{1}-\varepsilon_{2}-\varepsilon_{3}-\varepsilon_{4}\right) / 2$. The generator matrix $M(F_4)$ for $F_4$ is
$$
M({F_4})=\left(\begin{array}{l}
\alpha_{1} \\
\alpha_{2} \\
\alpha_{3} \\
\alpha_{4}
\end{array}\right)=\left(\begin{array}{cccc}
 & 1 & -1 &  \\
 &  & 1 & -1 \\
 &  & & 1 \\
\frac{1}{2} & -\frac{1}{2} & -\frac{1}{2} & -\frac{1}{2}
\end{array}\right)
$$
and its Gram matrix is
$$
G = M(F_4)\cdot{}^t M(F_4)=\left(\begin{array}{rrrr}
2 & -1 &  &  \\
-1 & 2 & -1 &  \\
 & -1 & 1 & -\frac{1}{2} \\
 &  & -\frac{1}{2} & 1
\end{array}\right)
$$
Also its theta function is (See Appendix 5.6 for a description of what a lattice theta function is.)
$$
\begin{aligned}
\Theta_{{F_{4}}}(q):=\sum_{x \in \Lambda_{F_{4}}} q^{(x, x)}=1+24 q+24 q^{2}+96 q^{3}+2 4 q^{4}+144 q^{5}+\ldots
\end{aligned}
$$
however $(x,x)$ to be Euclidean norm and $\Lambda_{F_{4}}$ to be the lattice generated by the basis of $F_4$.
\subsubsection{Root system $E_8$}
Let $\left\{\varepsilon_{1}, \varepsilon_{2}, \varepsilon_{3}, \varepsilon_{4},\varepsilon_{5},\varepsilon_{6},\varepsilon_{7},\varepsilon_{8}\right\}$ be the standard basis of
$\mathbb{R}^8$. I set the notation for $E_8$ as generated by the basis $\alpha_{i}:=\varepsilon_{i-1}-\varepsilon_{i-2}(1 \leq i \leq 6),\alpha_{7}:=\varepsilon_{1}+\varepsilon_{2} ,\alpha_{8}:=\frac{1}{2}\left(\varepsilon_{1}-\varepsilon_{2}-\varepsilon_{3}-\varepsilon_{4}-\varepsilon_{5}-\varepsilon_{6}-\varepsilon_{7}+\varepsilon_{8}\right)$
The generator matrix $M(E_8)$ for $E_8$ is
$$
M({E_8})=\left(\begin{array}{l}
\alpha_{1} \\
\alpha_{2} \\
\alpha_{3} \\
\alpha_{4} \\
\alpha_{5} \\
\alpha_{6} \\
\alpha_{7} \\
\alpha_{8}
\end{array}\right)=\left(\begin{array}{cccccccc}
1 & -1 &  &  &  &  &  &  \\
 & 1 & -1 &  &  &  &  &  \\
 &  & 1 & -1 &  &  &  &  \\
 &  &  & 1 & -1 &  &  &  \\
 &  &  &  & 1 & -1 &  &  \\
 &  &  &  &  & 1 & -1 &  \\
 &  &  &  &  & 1 & 1 &  \\
1 / 2 & -1 / 2 & -1 / 2 & -1 / 2 & -1 / 2 & -1 / 2 & -1 / 2 & -1 / 2
\end{array}\right)
$$
and its Gram matrix is
$$
G(E_8)=M(E_8)\cdot {}^t M(E_8)= \left(\begin{array}{rrrrrrrr}
2 & -1 &  &  &  &  &  &  \\
-1 & 2 & -1 &  &  &  &  &  \\
 & -1 & 2 & -1 &  &  &  &  \\
 &  & -1 & 2 & -1 &  &  &  \\
 &  &  & -1 & 2 & -1 & -1 &  \\
 &  &  &  & -1 & 2 &  &  \\
 &  &  &  & -1 &  & 2 & -1 \\
 &  &  &  &  &  & -1 & 2
\end{array}\right)
$$
Also its theta function is

$$
\begin{aligned}
\Theta_{E_{8}}(z) &=\sum_{x \in \Lambda_{E_8}}q^{(x,x)}\\
&=\frac{1}{2}\left(\theta_{2}(q)^{8}+\theta_{3}(q)^{8}+\theta_{4}(q)^8\right) \\
&=\theta_{2}(q^2)^{8}+14 \theta_{2}(q^2)^{4} \theta_{3}(q^2)^{4}+\theta_{3}(q^2)^{8}\\
&=1  +240q^2 + 2160q^4 + 6720q^6  + 17520q^8 + \cdots
\end{aligned}
$$
The theta function of $E_8$ is well known as Eisenstein series $E_4(z) = 1+240\sum_{n=1}^\infty \sigma_3(n) q^n$ in the theory of modular form.[1, P122] However,
$$
\begin{aligned}
\sigma_{k-1}(n)&=\sum_{d \mid n} d^{k-1}\\
\theta_{2}(q)&=\sum_{n=-\infty}^{\infty} q^{(m+1 / 2)^{2}} \\
\theta_{3}(q)&=\sum_{n=-\infty}^{\infty} q^{m^{2}} \\
\theta_{4}(q)&=\sum_{n=-\infty}^{\infty}(-q)^{m^{2}}
\end{aligned}
$$
Weyl group $W(F_4)$ and $W(E_8)$ are generated by reflections through the hyperplanes orthogonal to the generators(roots) of $F_4$ and $E_8$ respectively.In addition the lattices $F_4, E_8$ can be constructed using Quaternions.[See Appendix 5.4 and 5.5]

\section{Dual Quaternion and Weyl Group $W\left(F_{4}\right)$}
For Hurwitz quaternionic integers $\mathscr{H}$ described above, I studied its direct product $\mathscr{H}^2$, in other words the group ring $\mathbb{Z}^\prime[\mathbb{Q}^2_8]$ over $\mathbb{Z}^\prime = \mathbb{Z} \oplus (\mathbb{Z} + 1/2)\oplus (\mathbb{Z} + 1/4)$  of
$$
\mathbb{Q}^2_8 := <i,j,p,q|ip=pi,iq=qi,jp=pj>
$$
where $i,j$ and $p,q$ satisfy
$$
\begin{aligned}
&i^{2}=j^{2}=k^{2}=i j k=-1 \\
&i j=-j i=k, \quad j k=-k j=i, \quad k i=-i k=j
\end{aligned}
$$
$$
\begin{aligned}
&p^{2}=q^{2}=r^{2}=p q r=-1 \\
&p q=-q p=r, \quad q r=-r q=p, \quad r p=-p r=q
\end{aligned}
$$
$<a,b>$ denotes here the group generated by the multiplication of $a$ and $b$.\\\\
There is the permutation representation of $\mathbb{Q}^2_8$ generated by following.(See Appendix 5.1.)
$$
\begin{aligned}
&(1,3,15,13)(2,11,16,5)(4,7,14,9)(6,10,12,8) \\
&(1,12,15,6)(2,9,16,7)(3,10,13,8)(4,11,14,5) \\
&(1,8,15,10)(2,5,16,11)(3,6,13,12)(4,7,14,9) \\
&(1,13,15,3),(2,9,16,7),(4,5,14,11),(6,10,12,8)
\end{aligned}
$$
or
$$
\begin{aligned}
&(1,2,7,8)(3,4,5,6)\\
&(1,5,7,3)(2,4,8,6)\\
&(1,6,7,4)(2,3,8,5)\\
&(1,8,7,2)(3,4,5,6)
\end{aligned}
$$
\begin{figure}[h]
 \centering
 \includegraphics[keepaspectratio, scale=0.4]
      {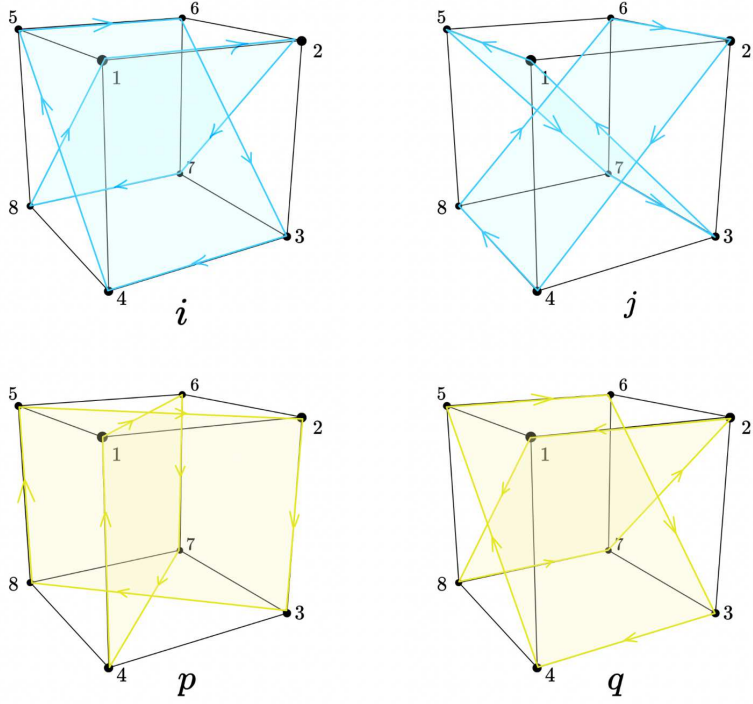}
 \caption{The permutations isomorphic to the generators of $\mathbb{Q}^2_8 $}
 \label{}
\end{figure}

For all elements $x$ in $\mathbb{Z}^\prime[\mathbb{Q}^2_8 ]$ can be written as
$$
x=\left(1, p, q, r\right)\left(\begin{array}{cccc}
a_{00} & a_{01} & a_{02} & a_{03} \\
a_{10} & a_{11} & a_{12} & a_{13} \\
a_{20} & a_{21} & a_{22} & a_{23} \\
a_{30} & a_{31} & a_{32} & a_{33}
\end{array}\right)\left(\begin{array}{l}
1 \\
i \\
j \\
k
\end{array}\right) \in \mathscr{H}^2
$$
where
$$
a_{ij} \in \mathbb{Z}^\prime =  \mathbb{Z} \oplus (\mathbb{Z} + 1/2) \oplus (\mathbb{Z} + 1/4)
$$
and (donate $\mathbb{Q} = <i,j>, \mathbb{Q}^\prime = <p,q>$, i.e. $\mathbb{Q}^2_8 \cong \mathbb{Q} \times \mathbb{Q}^\prime$)
$$
x \cdot x^\prime  = x^\prime  \cdot x \quad ({}^\forall x \in \mathbb{Q}, {}^\forall x^\prime \in \mathbb{Q}^\prime )
$$
The conjugation of $x \in \mathbb{Z}^{\prime}\left[\mathbb{Q}^2_8\right]$  follows the general one for group rings described above
$$
\bar{x}=\left(1, p, q, r\right)\left(\begin{array}{cccc}
a_{00} & -a_{01} & -a_{02} & -a_{03} \\
-a_{10} & a_{11} & a_{12} & a_{13} \\
-a_{20} & a_{21} & a_{22} & a_{23} \\
-a_{30} & a_{31} & a_{32} & a_{33}
\end{array}\right)\left(\begin{array}{l}
1 \\
i \\
j \\
k
\end{array}\right) \in \mathscr{H}^2
$$
Unlike the Quaternion, the norm $N(x)=x\bar{x}$ now dealing with is not always real number, but only in some special cases. One of them is Conway and Slone's $\mathcal{S} \cong 2 \cdot A_4$, but it also turns out that there is another finite group order $1152 = 2^{7} \cdot 3^{2}$.
\begin{itembox}{Lemma.1}
$\rho: \mathbb{Q} \rightarrow SO(4,\mathbb{R})$ is isomorphism.
  $$
  \rho: i \quad \longmapsto\left(\begin{array}{cccc} & 1 &  &  \\ -1 &  &  &  \\  &  &  & 1 \\  &  & -1 & \end{array}\right)
  $$
  $$
  \rho:j \quad \longmapsto\left(\begin{array}{cccc}
   &  & 1 &  \\
   &  &  & -1 \\
  -1 &  &  &  \\
   & 1 &  &
  \end{array}\right)
  $$
  $$
  \rho:k\quad \longmapsto\left(\begin{array}{cccc}
   &  &  & -1 \\
   &  & -1 &  \\
   & 1 &  &  \\
  1 &  &  &
  \end{array}\right)
  $$
\end{itembox}
(Proof)
Since $\mathbb{Q} \cong \mathrm{Q}_{8}=\left\langle a, b \mid a^{4}=e, a^{2}=b^{2}, b a=a^{-1} b\right\rangle$, To show that they are isomorphic, it is enough to show that the generators $i,j$ satisfy $i^{4}=E, \quad i^{2}=j^{\prime 2}, j i=i^{-1} j=-i j$. Calculating
$$
(i)\rho(j)\rho=\left(\begin{array}{llll}
 & 1 &  &  \\
-1 &  &  &  \\
 &  &  & 1 \\
 &  & -1 &
\end{array}\right)\left(\begin{array}{cccc}
 &  & 1 &  \\
 &  &  & -1 \\
-1 &  &  &  \\
 & 1 &  &
\end{array}\right)=\left(\begin{array}{cccc}
   &  &  & -1 \\
   &  & -1 &  \\
   & 1 &  &  \\
  1 &  &  &
  \end{array}\right)=(k)\rho
$$

$$
((i)\rho)^2=((j)\rho)^2=((k)\rho)^2=\left(\begin{array}{rrrr}
-1 &  &  &  \\
 & -1 &  &  \\
 &  & -1 &  \\
 &  &  & -1
\end{array}\right)
$$

$$
(j)\rho(i)\rho=\left(\begin{array}{rrrr}
 &  &  & 1 \\
 &  & 1 &  \\
 & -1 &  &  \\
-1 &  &  &
\end{array}\right)=-(i)\rho(j)\rho=-(k)\rho
$$
$$
((i)\rho)^{-1}(j)\rho=\left(\begin{array}{rrrr}
 & -1 &  &  \\
1 &  &  &  \\
 &  &  & -1 \\
 &  & 1 &
\end{array}\right)\left(\begin{array}{rrrr}
 &  & 1 &  \\
 &  &  & -1 \\
-1 &  &  &  \\
 & 1 &  &
\end{array}\right)=\left(\begin{array}{rrrr}
 &  &  & 1 \\
 &  & 1 &  \\
 & -1 &  &  \\
-1 &  &  &
\end{array}\right)
$$
This result shows that the above is true for the generators.(Q.E.D)
\begin{itembox}{Theorem.1}
The multiplicative group $G$ order $1152 = 2^{7} \cdot 3^{2}$ generated by
$$
\begin{aligned}
&e_{1}=\frac{1}{2}(1+i q+j r-k p) \\
&e_{2}=\frac{1}{2}(1-i p-j q-k r) \\
&e_{3}=\frac{1}{2}(1-i p-j r+k q) \\
&e_{4}=\frac{1}{2}(1-i q-j p+k r)
\end{aligned}
$$
is isomorphic to the exceptional Weyl Group$\  W\left(F_{4}\right)$.
$$
G \cong W\left(F_{4}\right)
$$
\end{itembox}
(Proof)\\
According to "ATLAS", the group$W(F_4)$ is defined by
\begin{equation}
W\left(F_{4}\right)=<a, b, c, d \mid a^{2}=b^{2}=c^{2}=d^{2}=(a b)^{3}=[a, c]=[a, d]=(b c)^{4}=[b, d]=(c d)^{3}=1>
\end{equation}
to show that the group $G$ is isomorphic to this group, I only need to check that the generating system of $G$ satisfies the same properties. Using a linear representation for simplicity. $\rho:i,j,k \in \mathbb{Q} \mapsto \boldsymbol{i},\boldsymbol{j},\boldsymbol{k}\in SO(4,\mathbb{R}) $ can be written as
$$
\begin{aligned}
&\boldsymbol{i} := (i) \rho=\left(\begin{array}{cccc}
 & 1 &  &  \\
-1 &  &  &  \\
 &  &  & 1 \\
 &  & -1 &
\end{array}\right) \\\\
&\boldsymbol{j} :=(j) \rho=\left(\begin{array}{cccc}
 &  & 1 &  \\
 &  &  & -1 \\
-1 &  &  &  \\
 & 1 &  &
\end{array}\right) \\\\
&\boldsymbol{k}:=(k) \rho=\left(\begin{array}{cccc}
 &  &  & -1 \\
 &  & -1 &  \\
 & 1 &  &  \\
1 &  &  &
\end{array}\right)
\end{aligned}
$$
Now determine $p, q$, and $r$ to be commutative with these. Let $(p)\rho$ as
$$
\boldsymbol{p}=(p)\rho=\left(\begin{array}{llll}
a_{11} & a_{12} & a_{13} & a_{14} \\
a_{21} & a_{22} & a_{23} & a_{24} \\
a_{31} & a_{32} & a_{33} & a_{34} \\
a_{41} & a_{42} & a_{43} & a_{44}
\end{array}\right)
$$
Since ${ }^{\forall} x \in \mathbb{Q},{ }^{\forall} x^{\prime} \in \mathbb{Q}^{\prime} $must be
 $x \cdot x^{\prime}=x^{\prime} \cdot x$, $\boldsymbol{p},\boldsymbol{i}$ and $\boldsymbol{j}$must satisfy
$\boldsymbol{p} \cdot \boldsymbol{i} =\boldsymbol{i} \cdot \boldsymbol{p}$ and $\boldsymbol{p}\cdot \boldsymbol{j} =\boldsymbol{i} \cdot \boldsymbol{p}$. Thus
$$
\begin{array}{c}
\boldsymbol{i}^{-1}\cdot  \boldsymbol{p}\cdot  \boldsymbol{i}=\left(\begin{array}{cccc}
a_{22} & -a_{21} & a_{24} & -a_{23} \\
-a_{12} & a_{11} & -a_{14} & a_{13} \\
a_{42} & -a_{41} & a_{44} & -a_{43} \\
-a_{32} & a_{31} & -a_{34} & a_{33}
\end{array}\right)=\boldsymbol{p} \\\\
\boldsymbol{j}^{-1}\cdot  \boldsymbol{p}\cdot  \boldsymbol{j}=\left(\begin{array}{cccc}
a_{33} & -a_{34} & -a_{31} & a_{32} \\
-a_{43} & a_{44} & a_{41} & -a_{42} \\
-a_{13} & a_{14} & a_{11} & -a_{12} \\
a_{23} & -a_{24} & -a_{21} & a_{22}
\end{array}\right)=\boldsymbol{p}
\end{array}
$$
Using appropriate $x,y,z,w$, $\boldsymbol{p}$ can be rewiriten as
$$
\boldsymbol{p}=\left(\begin{array}{cccc}
x & y & z & -w \\
-y & x & w & z \\
-z & -w & x & -y \\
w & -z & y & x
\end{array}\right)
$$
Also for $\boldsymbol{p}^2 = -E$.
$$
\begin{array}{l}
\left(\begin{array}{cccc}
x & y & -z & w \\
-y & x & w & -z \\
z & -w & x & -y \\
w & z & y & x
\end{array}\right)\left(\begin{array}{cccc}
x & y & -z & w \\
-y & x & w & -z \\
z & -w & x & -y \\
w & z & y & x
\end{array}\right)=\\
\left(\begin{array}{cccc}
x^{2}+w^{2}-\left(y^{2}+z^{2}\right), & 2(x y+z w), & 2(-x z+y w), & 2(x w-y z) \\
2(-x y+w z) & x^{2}-y^{2}-z^{2}-w^{2}, & 2(x w) & -2(y w+x z) \\
2 x z & -2 x w & x^{2}-y^{2}-z^{2}-w^{2} & 2(w z-x y) \\
2 x w & 2(x z) & 2 x y . & x^{2}-y^{2}-z^{2}+w^{2}
\end{array}\right)=-E
\end{array}
$$
where $E$ is the unit matrix of $4 \times 4$.To satisfy the above, $x = w = 0$.Rewriting
$$
\boldsymbol{p}^{2}=\left(\begin{array}{cccc}
-\left(y^{2}+z^{2}\right) &  &  & 2(-y z) \\
 & -\left(y^{2}+z^{2}\right) &  &  \\
 &  & -\left(y^{2}+z^{2}\right) &  \\
 &  &  & -\left(y^{2}+z^{2}\right)
\end{array}\right)
$$
This shows that either $y,z$ is one of $0,\pm 1$ respectively.From this we obtain $\boldsymbol{p},\boldsymbol{q},\boldsymbol{r}$ as follows.
\begin{equation}
\begin{array}{l}
\boldsymbol{p}=(p)\rho=\left(\begin{array}{cccc}
 & 1 &  &  \\
-1 &  &  &  \\
 &  &  & -1 \\
 &  & 1 &
\end{array}\right) \\\\
\boldsymbol{q}=(q)\rho=\left(\begin{array}{cccc}
 &  & 1 &  \\
 &  &  & 1 \\
-1 &  &  &  \\
 & -1 &  &
\end{array}\right) \\\\
\boldsymbol{r}=(r)\rho=\left(\begin{array}{cccc}
 &  &  & 1 \\
 &  & -1 &  \\
 & 1 &  &  \\
-1 &  &  &
\end{array}\right)
\end{array}
\end{equation}
Now we can write the generators as a matrix by $\rho$,
\begin{equation}
\begin{array}{l}
\boldsymbol{e}_1=(e_{1})\rho=\frac{1}{2}\left(\begin{array}{cccc}
1 & 1 & 1 & 1 \\
1 & 1 & -1 & -1 \\
1 & -1 & 1 & -1 \\
1 & -1 & -1 & 1
\end{array}\right) \\\\
\boldsymbol{e}_2=(e_{2})\rho=\left(\begin{array}{llll}
1 &  &  &  \\
 & 1 &  &  \\
 &  & 1 &  \\
 &  &  & -1
\end{array}\right) \\\\
\boldsymbol{e}_3=(e_{3})\rho=\left(\begin{array}{llll}
1 &  &  &  \\
 & 1 &  &  \\
 &  &  & 1 \\
 &  & 1 &
\end{array}\right) \\\\
\boldsymbol{e}_4=(e_{4})\rho=\left(\begin{array}{llll}
1 &  &  &  \\
 &  & 1 &  \\
 & 1 &  &  \\
 &  &  & 1
\end{array}\right)
\end{array}
\end{equation}
It can be shown by a simple calculation that these $\boldsymbol{e}_1,\boldsymbol{e}_2,\boldsymbol{e}_3,\boldsymbol{e}_4$ satisfy (1).Actually
\begin{equation}
\boldsymbol{e}_1^{2}=\boldsymbol{e}_2^{2}=\boldsymbol{e}_3^{2}=\boldsymbol{e}_4^{2}=\left(\boldsymbol{e}_1 \boldsymbol{e}_2\right)^{2}=\left[\boldsymbol{e}_1, \boldsymbol{e}_3\right]=\left[\boldsymbol{e}_2, \boldsymbol{e}_4\right]=\left[\boldsymbol{e}_1, \boldsymbol{e}_4\right]=\left(\boldsymbol{e}_2 \boldsymbol{e}_3\right)^{4}=\left(\boldsymbol{e}_3 \boldsymbol{e}_4\right)^{3}=1
\end{equation}
However, the ecommutator is defined to be $\left[\boldsymbol{e}_{i}, \boldsymbol{e}_{j}\right]=\boldsymbol{e}_{i} \cdot \boldsymbol{e}_{j} \cdot \boldsymbol{e}_{i}^{-1} \cdot \boldsymbol{e}_{j}^{-1}$.\\
(Q.E.D.)
\section{The $16$ dimensional lattice $\Lambda_{\mathbb{Q}^2_8}$ and Barns-Wall lattice}
According to J.H.Conway and Neil J. A. Sloane[1], the set of norm $1$ of quaternions generated by the multiplication with follwings  forms the group isomorphic to $2 \cdot A_4$ order $24$.
$$
<i,j,k,\frac{1+i+j+k}{2}>
$$
And by the inverse mapping $\phi$ which is the bijection
$$
\begin{array}{llllll}
\phi: &(\mathbb{Z}\oplus \mathbb{Z}+1/2)^4& \rightarrow& \left \langle i,j,k,\frac{1+i+j+k}{2} \right \rangle& \\
&(c_0,c_1,c_2,c_3)& \mapsto & c_0 + c_1 i + c_2 j + c_3 k&
\end{array}
$$
we obtain $F_4$ lattice[See Appendix $5.4$] whose theta series is
$$
\Theta_{F_4}(q)=2 E_{2}(q^2)-E_{2}(q)=1+24 q+24 q^{2}+96 q^{3}+24 q^{4}+144 q^{5}+\ldots
$$
where $E_2$ is Eisenstein series
$$
E_{2}(q)=1-24 \sum_{n} \sigma_{1}(n) q^{n}
$$

The group $G \cong W(F_4)$ is subring of the group ring $\mathbb{Z}^{\prime}\left[\mathbb{Q}^2_8\right]$ of norm 1. So let donate $\varphi$ the bijection $\varphi: (\mathbb{Z} \oplus \mathbb{Z}+1 / 2\oplus \mathbb{Z}+1 / 4)^{16} {\longrightarrow} \mathbb{Z}^{\prime}\left[\mathbb{Q}^2_8\right]$
and donate $N$ norm $N: x \rightarrow x \bar{x}$, then $G$ can be also regarded as the subset of $\text{Ker}(N)$.
\\\\
In terms of the composition of the following maps,
$$
\begin{array}{llllllll}
&(\mathbb{Z} \oplus \mathbb{Z}+1 / 2\oplus \mathbb{Z}+1 / 4)^{16}& \stackrel{\varphi}{\rightarrow}&\quad \quad \quad  \mathbb{Z}^{\prime}\left[\mathbb{Q}^2_8\right]& \stackrel{N}{\rightarrow}&\quad \quad \mathbb{Z}^{\prime}\left[\mathbb{Q}^2_8\right]&\\
&\quad \quad \quad \rotatebox{90}{$\in$}& &\quad \quad \quad \quad \quad \rotatebox{90}{$\in$}& &\quad \quad \quad \quad  \rotatebox{90}{$\in$}&\\
&\left(x_{1}, x_{2}, \ldots x_{16}\right) &\longmapsto & a_{00}+a_{01} i+\cdots+a_{33} k r& \mapsto& a_{00}^{\prime}+a_{01}^{\prime}i \cdots a_{33}^{\prime}k r&
\end{array}
$$

\begin{itembox}{Proposition.1}
The freemodule generated by the following elements of $G$ is the lattice of rank $16$
$$
\begin{aligned}
&\frac{1}{2}(1+i+j+k) \\
&\frac{1}{2}(1+i+j-k) \\
&\frac{1}{2}(1+i)(1+q) \\
&\frac{1}{2}(1+i)(1+r) \\
&\frac{1}{2}(1+j)(1+q) \\
&\frac{1}{2}(1+j)(1+r) \\
&\frac{1}{2}(1+i q+j r-p k) \\
&\frac{1}{2}(1-i r+j q-k p) \\
&\frac{1}{2}(1-p+q-r) \\
&\frac{1}{2}(1-p-q+r) \\
&\frac{1}{2}(1-k)(1+q) \\
&\frac{1}{2}(1-k)(1+r) \\
&\frac{1}{2}(1+i r-j q-k p) \\
&\frac{1}{2}(1+j)(1-r)\\
&\frac{1}{2}(1+j)(1-p)\\
&\frac{1}{4}(1 + i + j + k + p + ip + \cdots + kr)
\end{aligned}
$$
\end{itembox}
(Proof)\\
The above elements of the $G$ can be written by bijection $\varphi^{-1}: G \rightarrow (\mathbb{Z} \oplus \mathbb{Z} + 1/2\oplus \mathbb{Z} + 1/4)^{16}$ as
$$
\varphi^{-1}:(1, p, q, r)\left(\begin{array}{llll}
a_{00} & a_{01} & a_{02} & a_{03} \\
a_{10} & a_{11} & a_{12} & a_{13} \\
a_{20} & a_{21} & a_{22} & a_{23} \\
a_{30} & a_{31} & a_{32} & a_{33}
\end{array}\right)\left(\begin{array}{l}
1 \\
i \\
j \\
k
\end{array}\right)\mapsto(a_{00},a_{01},a_{02},a_{03},a_{10},a_{11},a_{12},a_{13},a_{20},...,a_{33})
$$
using the map, all the generator of the module are shown as
$$
\left(\begin{aligned}
&\frac{1}{2}(1+i+j+k) \\
&\frac{1}{2}(1+i+j-k) \\
&\frac{1}{2}(1+i)(1+q) \\
&\frac{1}{2}(1+i)(1+r) \\
&\frac{1}{2}(1+j)(1+q) \\
&\frac{1}{2}(1+j)(1+r) \\
&\frac{1}{2}(1+i q+j r-p k) \\
&\frac{1}{2}(1-i r+j q-k p) \\
&\frac{1}{2}(1-p+q-r) \\
&\frac{1}{2}(1-p-q+r) \\
&\frac{1}{2}(1-k)(1+q) \\
&\frac{1}{2}(1-k)(1+r) \\
&\frac{1}{2}(1+i r-j q-k p) \\
&\frac{1}{2}(1+j)(1-r) \\
&\frac{1}{2}(1+j)(1-p) \\
&\frac{1}{4}(1+i+\cdots+k r)
\end{aligned}\right) \mapsto N=
\frac{1}{4}
\left(\begin{array}{rrrrrrrrrrrrrrrr}
2 & 2 & 2 & 2 & 0 & 0 & 0 & 0 & 0 & 0 & 0 & 0 & 0 & 0 & 0 & 0 \\
2 & 2 & 2 & -2 & 0 & 0 & 0 & 0 & 0 & 0 & 0 & 0 & 0 & 0 & 0 & 0 \\
2 & 2 & 0 & 0 & 0 & 0 & 0 & 0 & 2 & 2 & 0 & 0 & 0 & 0 & 0 & 0 \\
2 & 2 & 0 & 0 & 0 & 0 & 0 & 0 & 0 & 0 & 0 & 0 & 2 & 2 & 0 & 0 \\
2 & 0 & 2 & 0 & 0 & 0 & 0 & 0 & 2 & 0 & 2 & 0 & 0 & 0 & 0 & 0 \\
2 & 0 & 2 & 0 & 0 & 0 & 0 & 0 & 0 & 0 & 0 & 0 & 2 & 0 & 2 & 0 \\
2 & 0 & 0 & 0 & 0 & 0 & 0 & -2 & 0 & 2 & 0 & 0 & 0 & 0 & 2 & 0 \\
2 & 0 & 0 & 0 & 0 & 0 & 0 & -2 & 0 & 0 & 2 & 0 & 0 & -2 & 0 & 0 \\
2 & 0 & 0 & 0 & -2 & 0 & 0 & 0 & 2 & 0 & 0 & 0 & -2 & 0 & 0 & 0 \\
2 & 0 & 0 & 0 & -2 & 0 & 0 & 0 & -2 & 0 & 0 & 0 & 2 & 0 & 0 & 0 \\
2 & 0 & 0 & -2 & 0 & 0 & 0 & 0 & 2 & 0 & 0 & -2 & 0 & 0 & 0 & 0 \\
2 & 0 & 0 & -2 & 0 & 0 & 0 & 0 & 0 & 0 & 0 & 0 & 2 & 0 & 0 & -2 \\
2 & 0 & 0 & 0 & 0 & 0 & 0 & -2 & 0 & 0 & -2 & 0 & 0 & 2 & 0 & 0 \\
2 & 0 & 2 & 0 & 0 & 0 & 0 & 0 & 0 & 0 & 0 & 0 & -2 & 0 & -2 & 0 \\
2 & 0 & 2 & 0 & -2 & 0 & -2 & 0 & 0 & 0 & 0 & 0 & 0 & 0 & 0 & 0 \\
1 & 1 & 1 & 1 & 1 & 1 & 1 & 1 & 1 & 1 & 1 & 1 & 1 & 1 & 1 & 1
\end{array}\right)
$$
to show the rank of the freemodule, I can confirm it by using line elementary row transformation of $N$. Actually $N$ can be transformed as
$$
M = \frac{1}{4}\left(\begin{array}{llllllllllllllll}
1 & 1 & 1 & 1 & 1 & 1 & 1 & 1 & 1 & 1 & 1 & 1 & 1 & 1 & 1 & 1 \\
 & 2 &  &  &  &  &  & 2 &  &  &  & 2 &  & 2 &  &  \\
 &  & 2 &  &  &  &  & 2 &  &  &  & 2 &  &  & 2 &  \\
 &  &  & 2 &  &  &  & 2 &  &  &  & 2 &  &  &  & 2 \\
 &  &  &  & 2 &  &  & 2 &  &  &  &  &  & 2 & 2 &  \\
 &  &  &  &  & 2 &  & 2 &  &  &  &  &  & 2 &  & 2 \\
 &  &  &  &  &  & 2 & 2 &  &  &  &  &  &  & 2 & 2 \\
 &  &  &  &  &  &  & 4 &  &  &  &  &  &  &  &  \\
 &  &  &  &  &  &  &  & 2 &  &  & 2 &  & 2 & 2 &  \\
 &  &  &  &  &  &  &  &  & 2 &  & 2 &  & 2 &  & 2 \\
 &  &  &  &  &  &  &  &  &  & 2 & 2 &  &  & 2 & 2 \\
 &  &  &  &  &  &  &  &  &  &  & 4 &  &  &  &  \\
 &  &  &  &  &  &  &  &  &  &  &  & 2 & 2 & 2 & 2 \\
 &  &  &  &  &  &  &  &  &  &  &  &  & 4 &  &  \\
 &  &  &  &  &  &  &  &  &  &  &  &  &  & 4 &  \\
 &  &  &  &  &  &  &  &  &  &  &  &  &  &  & 4
\end{array}\right)
$$
Therefore the rank of the freemodule is $16$.(Q.E.D.)\\\\

\begin{itembox}{Proposition.2}
  Let $\Lambda_{\mathbb{Q}^2_8}$ be the lattice generated by the above generator.Thus
  $$
  \Lambda_{\mathbb{Q}^2_8}=\text{Span}_{\mathbb{Z}}\frac{1}{4}
\left(\begin{array}{llllllllllllllll}
1 & 1 & 1 & 1 & 1 & 1 & 1 & 1 & 1 & 1 & 1 & 1 & 1 & 1 & 1 & 1 \\
 & 2 &  &  &  &  &  & 2 &  &  &  & 2 &  & 2 &  &  \\
 &  & 2 &  &  &  &  & 2 &  &  &  & 2 &  &  & 2 &  \\
 &  &  & 2 &  &  &  & 2 &  &  &  & 2 &  &  &  & 2 \\
 &  &  &  & 2 &  &  & 2 &  &  &  &  &  & 2 & 2 &  \\
 &  &  &  &  & 2 &  & 2 &  &  &  &  &  & 2 &  & 2 \\
 &  &  &  &  &  & 2 & 2 &  &  &  &  &  &  & 2 & 2 \\
 &  &  &  &  &  &  & 4 &  &  &  &  &  &  &  &  \\
 &  &  &  &  &  &  &  & 2 &  &  & 2 &  & 2 & 2 &  \\
 &  &  &  &  &  &  &  &  & 2 &  & 2 &  & 2 &  & 2 \\
 &  &  &  &  &  &  &  &  &  & 2 & 2 &  &  & 2 & 2 \\
 &  &  &  &  &  &  &  &  &  &  & 4 &  &  &  &  \\
 &  &  &  &  &  &  &  &  &  &  &  & 2 & 2 & 2 & 2 \\
 &  &  &  &  &  &  &  &  &  &  &  &  & 4 &  &  \\
 &  &  &  &  &  &  &  &  &  &  &  &  &  & 4 &  \\
 &  &  &  &  &  &  &  &  &  &  &  &  &  &  & 4
\end{array}\right)
$$
then it satisfies
$$
\Lambda_{16} \subset \Lambda_{\text{Barns-Wall}}
$$
where $\Lambda_{\text{Barns-Wall}}$ generated below is called the Barnes-Wall lattice.[1.P131]
  $$
  \Lambda_{\text{Barns-Wall}}=\text{Span}_{\mathbb{Z}}\frac{1}{4}\left(\begin{array}{llllllllllllllll}
1 &  &  &  &  & 1 &  & 1 &  &  & 1 & 1 &  & 1 & 1 & 1 \\
 & 1 &  &  &  & 1 & 1 & 1 & 1 &  & 1 &  & 1 & 1 &  &  \\
 &  & 1 &  &  &  & 1 & 1 & 1 & 1 &  & 1 &  & 1 & 1 &  \\
 &  &  & 1 &  & 1 &  &  & 1 & 1 &  & 1 & 1 & 1 &  & 1 \\
 &  &  &  & 1 &  & 1 &  &  & 1 & 1 &  & 1 & 1 & 1 & 1 \\
 &  &  &  &  & 2 &  &  &  &  &  &  &  &  &  & 2 \\
 &  &  &  &  &  & 2 &  &  &  &  &  &  &  &  & 2 \\
 &  &  &  &  &  &  & 2 &  &  &  &  &  &  &  & 2 \\
 &  &  &  &  &  &  &  & 2 &  &  &  &  &  &  & 2 \\
 &  &  &  &  &  &  &  &  & 2 &  &  &  &  &  & 2 \\
 &  &  &  &  &  &  &  &  &  & 2 &  &  &  &  & 2 \\
 &  &  &  &  &  &  &  &  &  &  & 2 &  &  &  & 2 \\
 &  &  &  &  &  &  &  &  &  &  &  & 2 &  &  & 2 \\
 &  &  &  &  &  &  &  &  &  &  &  &  & 2 &  & 2 \\
 &  &  &  &  &  &  &  &  &  &  &  &  &  & 2 & 2 \\
 &  &  &  &  &  &  &  &  &  &  &  &  &  &  & 4
\end{array}\right)
  $$
\end{itembox}
(Proof)\\
For all elements $x \in \Lambda_{\mathbb{Q}^2_8}$ can be written by the following linear combination.
$$
x = (c_1,c_2,...,c_{16}) \cdot M
$$
where $(c_1,c_2,...,c_{16}) \in \mathbb{Z}^{16}$ and $M$ is the generator matrix of $\Lambda_{\mathbb{Q}^2_8}$. Let $M^\prime$ be the generator matrix of $\Lambda_{\text{Barns-Wall}}$ and if $\Lambda_{16} \subset \Lambda_{\text{Barns-Wall}}$ then for all $c \in \mathbb{Z}^{16}$, there exists $(a_{1},a_{2},...,a_{16})$
such that
$$
(c_1,c_2,...,c_{16})M = (a_1,a_2,...,a_{16})M^\prime
$$
in other words, the following must always be met if $\Lambda_{\mathbb{Q}^2_8} \subset \Lambda_{\text{Barns-Wall}}$
\begin{equation}
(c_1,c_2,...,c_{16})\cdot M \cdot (M^\prime)^{-1} \in \mathbb{Z}^{16}
\end{equation}
By actually performing the calculation of $M\cdot (M^\prime)^{-1}$
$$
M\cdot (M^\prime)^{-1}=\left(\begin{array}{rrrrrrrrrrrrrrrr}
1 & 1 & 1 & 1 & 1 & -1 & -1 & -1 & -1 & -1 & -1 & -1 & -1 & -2 & -1 & 5 \\
 & 2 &  &  &  & -1 & -1 &  & -1 &  & -1 & 1 & -1 &  &  & 2 \\
 &  & 2 &  &  &  & -1 &  & -1 & -1 &  &  &  & -1 &  & 2 \\
 &  &  & 2 &  & -1 &  & 1 & -1 & -1 &  &  & -1 & -1 &  & 2 \\
 &  &  &  & 2 &  & -1 & 1 &  & -1 & -1 &  & -1 &  &  & 1 \\
 &  &  &  &  & 1 &  & 1 &  &  &  &  &  & 1 &  & -1 \\
 &  &  &  &  &  & 1 & 1 &  &  &  &  &  &  & 1 & -1 \\
 &  &  &  &  &  &  & 2 &  &  &  &  &  &  &  & -1 \\
 &  &  &  &  &  &  &  & 1 &  &  & 1 &  & 1 & 1 & -2 \\
 &  &  &  &  &  &  &  &  & 1 &  & 1 &  & 1 &  & -1 \\
 &  &  &  &  &  &  &  &  &  & 1 & 1 &  &  & 1 & -1 \\
 &  &  &  &  &  &  &  &  &  &  & 2 &  &  &  & -1 \\
 &  &  &  &  &  &  &  &  &  &  &  & 1 & 1 & 1 & -1 \\
 &  &  &  &  &  &  &  &  &  &  &  &  & 2 &  & -1 \\
 &  &  &  &  &  &  &  &  &  &  &  &  &  & 2 & -1 \\
 &  &  &  &  &  &  &  &  &  &  &  &  &  &  & 1
\end{array}\right)
$$
Obviously, equation (5) is satisfied since $(c_1,c_2,...,c_{16}) \in \mathbb{Z}^{16}$.Therefore $\Lambda_{\mathbb{Q}^2_8} \subset \Lambda_{\text{Barns-Wall}}$\\
(Q.E.D)\\\\
By the above, the following is shown at the same time.
$$
\Lambda_{\text{Barns-Wall}} \not \subset \Lambda_{\mathbb{Q}^2_8}
$$
According to Conway and Slone,As the theta series shows the lattice $\Lambda_{\text{Barns-Wall}}$ is known to be a 16-dimensional dense packing(however this is considered controversial.) and its kissing number is 4320.
$$
\begin{aligned}
\Theta_{\Lambda_{\text{Barnes-Wall}}}(z) &=1 / 2\left\{\theta_{2}(q^2)^{16}+\theta_{3}(q^2)^{16}+\theta_{4}(q^2)^{16}+30 \theta_{2}(q^2)^{8} \theta_{3}(q^2)^{8}\right\} \\
&=1+4320 q^{4}+61440 q^{6}+\cdots .
\end{aligned}
$$
It is also mentioned that $\Lambda_{\text{Barns-Wall}}$ may constructed from the Leech lattice $\Lambda_{24}$ which is invariant under the action of the Conwy group.\\\\
The Gram matrix of $M$ which is the generator matrix of $\Lambda_{\mathbb{Q}^2_8}$ is following by the way
\begin{equation}
G := M\cdot {}^tM=\frac{1}{4}\left(\begin{array}{llllllllllllllll}
4 & 2 & 2 & 2 & 2 & 2 & 2 & 1 & 2 & 2 & 2 & 1 & 2 & 1 & 1 & 1 \\
2 & 4 & 2 & 2 & 2 & 2 & 1 & 2 & 2 & 2 & 1 & 2 & 1 & 2 &  &  \\
2 & 2 & 4 & 2 & 2 & 1 & 2 & 2 & 2 & 1 & 2 & 2 & 1 &  & 2 &  \\
2 & 2 & 2 & 4 & 1 & 2 & 2 & 2 & 1 & 2 & 2 & 2 & 1 &  &  & 2 \\
2 & 2 & 2 & 1 & 4 & 2 & 2 & 2 & 2 & 1 & 1 &  & 2 & 2 & 2 &  \\
2 & 2 & 1 & 2 & 2 & 4 & 2 & 2 & 1 & 2 & 1 &  & 2 & 2 &  & 2 \\
2 & 1 & 2 & 2 & 2 & 2 & 4 & 2 & 1 & 1 & 2 &  & 2 &  & 2 & 2 \\
1 & 2 & 2 & 2 & 2 & 2 & 2 & 4 &  &  &  &  &  &  &  &  \\
2 & 2 & 2 & 1 & 2 & 1 & 1 &  & 4 & 2 & 2 & 2 & 2 & 2 & 2 &  \\
2 & 2 & 1 & 2 & 1 & 2 & 1 &  & 2 & 4 & 2 & 2 & 2 & 2 &  & 2 \\
2 & 1 & 2 & 2 & 1 & 1 & 2 &  & 2 & 2 & 4 & 2 & 2 &  & 2 & 2 \\
1 & 2 & 2 & 2 &  &  &  &  & 2 & 2 & 2 & 4 &  &  &  &  \\
2 & 1 & 1 & 1 & 2 & 2 & 2 &  & 2 & 2 & 2 &  & 4 & 2 & 2 & 2 \\
1 & 2 &  &  & 2 & 2 &  &  & 2 & 2 &  &  & 2 & 4 &  &  \\
1 &  & 2 &  & 2 &  & 2 &  & 2 &  & 2 &  & 2 &  & 4 &  \\
1 &  &  & 2 &  & 2 & 2 &  &  & 2 & 2 &  & 2 &  &  & 4
\end{array}\right)
\end{equation}
For the density of the Barns-wall lattice and the lattice $\Lambda_{\mathbb{Q}^2_8}$ is as follows.
\begin{equation}
\frac{\det(M)}{\det(M^\prime)}=256=(\sqrt{2})^{16}
\end{equation}

\begin{itembox}{Proposition.3}

For any element of $\Lambda^\prime_{\mathbb{Q}^2_8}=\{x \in \Lambda_{\mathbb{Q}^2_8}| 4 \cdot x\}$, the components are either all even or all odd.
\end{itembox}
(Proof)\\
The generator matrix $M_4$ of $\Lambda^\prime_{\mathbb{Q}^2_8}$ is defined to be
$$
M_4 = \left(\begin{array}{c}
e_{1} \\
e_{2} \\
e_{3} \\
e_{4} \\
\vdots \\
e_{16}
\end{array}\right) =
\left(\begin{array}{llllllllllllllll}
1 & 1 & 1 & 1 & 1 & 1 & 1 & 1 & 1 & 1 & 1 & 1 & 1 & 1 & 1 & 1 \\
 & 2 &  &  &  &  &  & 2 &  &  &  & 2 &  & 2 &  &  \\
 &  & 2 &  &  &  &  & 2 &  &  &  & 2 &  &  & 2 &  \\
 &  &  & 2 &  &  &  & 2 &  &  &  & 2 &  &  &  & 2 \\
 &  &  &  & 2 &  &  & 2 &  &  &  &  &  & 2 & 2 &  \\
 &  &  &  &  & 2 &  & 2 &  &  &  &  &  & 2 &  & 2 \\
 &  &  &  &  &  & 2 & 2 &  &  &  &  &  &  & 2 & 2 \\
 &  &  &  &  &  &  & 4 &  &  &  &  &  &  &  &  \\
 &  &  &  &  &  &  &  & 2 &  &  & 2 &  & 2 & 2 &  \\
 &  &  &  &  &  &  &  &  & 2 &  & 2 &  & 2 &  & 2 \\
 &  &  &  &  &  &  &  &  &  & 2 & 2 &  &  & 2 & 2 \\
 &  &  &  &  &  &  &  &  &  &  & 4 &  &  &  &  \\
 &  &  &  &  &  &  &  &  &  &  &  & 2 & 2 & 2 & 2 \\
 &  &  &  &  &  &  &  &  &  &  &  &  & 4 &  &  \\
 &  &  &  &  &  &  &  &  &  &  &  &  &  & 4 &  \\
 &  &  &  &  &  &  &  &  &  &  &  &  &  &  & 4
\end{array}\right)
$$
We can find that the generators of the lattice $e_1 \in (2\mathbb{Z}+1)^{16}$ and others $e_i \in (2 \mathbb{Z})^{16}\ (i=1,2,...,16)$.Therefore, all of the elements generated by their addition are also even or odd.\\
(Q.E.D)\\\\
From now, $\Lambda^\prime_{\mathbb{Q}^2_8}$ will be simply written as $\Lambda$.
\begin{itembox}{Proposition.4}
The following homomorphism exist.
$$
f: \Lambda \rightarrow 2\Lambda
$$
where
$$
2 \Lambda :=\left\{x \in \Lambda \mid \sum_{i=1}^{16} c_{i} e_{i},\quad c_{1} \in 2 \mathbb{Z},\quad c_{i} \in \mathbb{Z}\quad(i=2,3, \ldots, 16)\right\}
$$
\end{itembox}
(Proof)\\
All the elements $x \in \Lambda$ can be written as the linear combination of the basis $\sum_{i=1}^{16} c_i x_i \quad (c_i \in \mathbb{Z})$, and whether the $x$ is even or odd is determined by the evenness of $c_1$. By defining the homomorphism as follows so the map to be $f: \Lambda \rightarrow 2\Lambda$.
$$
\begin{array}{lllll}
  f : &\mathbb{Z}& \rightarrow &2\mathbb{Z}&\\
      &c_1 & \mapsto & c^\prime_1&
\end{array}
$$
($\mathbb{Z}/\text{Ker}f = \mathbb{Z}/(2\mathbb{Z}) \cong \mathbb{F}_2$.) In other words
$$
\begin{array}{lllll}
  f : &\mathbb{Z}^{16}& \rightarrow &(2\mathbb{Z})\oplus\mathbb{Z}^{15}&\\
      &(c_1,c_2,...,c_{16}) & \mapsto & (2c_1,c_2,...,c_{16}) &
\end{array}
$$
it is easily to show the map is homomrphic. For any $x,y \in \Lambda$ can be written as
$$
\begin{aligned}
  x &= c_1 e_1 + \sum_{i=2}^{16}c_i e_i\\
  y &= a_1 e_1 + \sum_{i=2}^{16}a_i e_i\\
\end{aligned}
$$
where $c,a \in \mathbb{Z}$. And from the definition
$$
\begin{aligned}
  f(x) &= 2c_1 e_1 + \sum_{i=2}^{16}c_i e_i\\
  f(y) &= 2a_1 e_1 + \sum_{i=2}^{16}a_i e_i\\
\end{aligned}
$$
on the other hand
$$
f(x + y) = 2(c_1+a_1) e_1 + \sum_{i=2}^{16}(c_i+a_1) e_i\\
$$
Thus
$$
f(x) + f(y) = f(x + y)
$$
(Q.E.D)
\begin{itembox}{Proposition.5}
For $x \in \Lambda$, the map $g$ follwing is the isomorphism.
$$
g: x \mapsto x/2
$$
\end{itembox}
(Proof)\\
It is obvious that $g$ is the bijection. For any $x,y \in \Lambda$,
$$
g(x) + g(y) = \frac{x}{2} + \frac{y}{2}
$$
on the other and
$$
g(x + y)  =\frac{x+y}{2}
$$
Thus
$$
g(x) + g(y) = g(x + y)
$$
hence $g$ is the homomorphism. Therefore $g$ is the isomorphism.\\
(Q.E.D)\\\\
Also for $x \in \Lambda$, donate $h$ to be the homomorphism
$$
h: x \mapsto x \quad (\text{mod}\ 2)
$$
then

\begin{itembox}{Proposition.6}
The composition $\pi=h\circ g \circ f$ is the homomorphism to the code (vector space of degree $16$ and dimension $11$ over Finite Field $\mathbb{F}_2$) $\mathcal{C}_{16}$ generated below from $\Lambda$.
$$
\mathcal{C}_{16}=\text{Span}_{\mathbb{F}_2}
\left(\begin{array}{llllllllllllllll}
1 & 0 & 0 & 0 & 0 & 0 & 0 & 1 & 0 & 0 & 0 & 1 & 0 & 1 & 1 & 1 \\
0 & 1 & 0 & 0 & 0 & 0 & 0 & 1 & 0 & 0 & 0 & 1 & 0 & 1 & 0 & 0 \\
0 & 0 & 1 & 0 & 0 & 0 & 0 & 1 & 0 & 0 & 0 & 1 & 0 & 0 & 1 & 0 \\
0 & 0 & 0 & 1 & 0 & 0 & 0 & 1 & 0 & 0 & 0 & 1 & 0 & 0 & 0 & 1 \\
0 & 0 & 0 & 0 & 1 & 0 & 0 & 1 & 0 & 0 & 0 & 0 & 0 & 1 & 1 & 0 \\
0 & 0 & 0 & 0 & 0 & 1 & 0 & 1 & 0 & 0 & 0 & 0 & 0 & 1 & 0 & 1 \\
0 & 0 & 0 & 0 & 0 & 0 & 1 & 1 & 0 & 0 & 0 & 0 & 0 & 0 & 1 & 1 \\
0 & 0 & 0 & 0 & 0 & 0 & 0 & 0 & 1 & 0 & 0 & 1 & 0 & 1 & 1 & 0 \\
0 & 0 & 0 & 0 & 0 & 0 & 0 & 0 & 0 & 1 & 0 & 1 & 0 & 1 & 0 & 1 \\
0 & 0 & 0 & 0 & 0 & 0 & 0 & 0 & 0 & 0 & 1 & 1 & 0 & 0 & 1 & 1 \\
0 & 0 & 0 & 0 & 0 & 0 & 0 & 0 & 0 & 0 & 0 & 0 & 1 & 1 & 1 & 1
\end{array}\right)
$$
\end{itembox}
(Proof)\\
From the Proposition.3, the lattice $\Lambda$ can be decomposed directly  $2\Lambda \oplus (2\Lambda+1)$, where
$$
\begin{aligned}
2 \Lambda &:=\left\{x \in \Lambda \mid \sum_{i=1}^{16} c_{i} e_{i},\quad c_{1} \in 2 \mathbb{Z},\quad c_{i} \in \mathbb{Z}\quad(i=2,3, \ldots, 16)\right\} \\
2 \Lambda+1 &:=\left\{x \in \Lambda \mid \sum_{i=1}^{16} c_{i} e_{i},\quad c_{1} \in (2 \mathbb{Z}+1),\quad c_{i} \in \mathbb{Z}\quad(i=2,3, \ldots, 16)\right\}
\end{aligned}
$$
and as I mentioned in Proposition.4 and 5 there are the homomorphism $f: \Lambda \rightarrow 2\Lambda$, and the isomorphism $g$, the homomorphism as well. Donate the images produced by these mappings to be as follows.
$$
\begin{array}{lllllllll}
&\mathbb{Z}^{16}&\stackrel{f}{\longrightarrow}&(2\mathbb{Z})^{16}& \stackrel{g}{\longrightarrow} &\mathbb{Z}^{16}&\stackrel{h}{\longrightarrow}&(\mathbb{F}_2)^{16}&\\
&\cup& &\cup& & \cup & & \cup &\\
&\Lambda & \stackrel{\text{Hom}}\longrightarrow & 2\Lambda & \stackrel{\text{Iso}}\longrightarrow & \Lambda^\prime & \stackrel{\text{Hom}}\longrightarrow & \Lambda(\mathbb{F}_2)&\\
\end{array}
$$
where $\mathbb{F}_2$ is the finite field.(The operation rules for this field is $1+1 = 0, 1+0 = 1, 0+0=0$)\\\\
For example the element $e_{13}$ will be the following by composition of the maps $\pi$.
$$
(0,0,0,...,0,2,2,2,2)\mapsto (0,0,0,...,0,2,2,2,2) \mapsto (0,0,0,...,0,1,1,1,1) \mapsto (0,0,0,...,0,1,1,1,1)
$$
Also $e_{16}$ and $e_1$ becomes
$$
\begin{aligned}
&(0,0,0,...,0,0,0,0,4)\mapsto (0,0,0,...,0,0,0,0,4) \mapsto (0,0,0,...,0,0,0,0,2) \mapsto (0,0,0,...,0,0,0,0,0)\\
&(1,1,1,...,1,1,1,1,1)\mapsto(2,2,2,...,2,2,2,2,2) \mapsto (1,1,1,...,1,1,1,1,1) \mapsto (1,1,1,...,1,1,1,1,1)
\end{aligned}
$$
So the generators become as follows by the homomorphism $\pi$.
$$
\pi: \left(\begin{array}{c}
e_{1}\\
e_{2} \\
e_{3} \\
e_{4} \\
e_{5} \\
\vdots \\
e_{16}
\end{array}\right)\mapsto
\left(\begin{array}{llllllllllllllll}
1 & 1 & 1 & 1 & 1 & 1 & 1 & 1 & 1 & 1 & 1 & 1 & 1 & 1 & 1 & 1\\
0 & 1 &  &  &  &  &  & 1 &  &  &  & 1 &  & 1 &  &  \\
\vdots &  & 1 &  &  &  &  & 1 &  &  &  & 1 &  &  & 1 &  \\
 &  &  & 1 &  &  &  & 1 &  &  &  & 1 &  &  &  & 1 \\
 &  &  &  & 1 &  &  & 1 &  &  &  &  &  & 1 & 1 &  \\
 &  &  &  &  & 1 &  & 1 &  &  &  &  &  & 1 &  & 1 \\
 \vdots &  &  &  &  &  & 1 & 1 &  &  &  &  &  &  & 1 & 1 \\
0 & \cdots &  &  &  &  &  &  &  &  &  &  &  &  &\cdots  & 0 \\
 \vdots &  &  &  &  &  &  &  & 1 &  &  & 1 &  & 1 & 1 &  \\
 &  &  &  &  &  &  &  &  & 1 &  & 1 &  & 1 &  & 1 \\
 \vdots &  &  &  &  &  &  &  &  &  & 1 & 1 &  &  & 1 & 1 \\
0 & \cdots &  &  &  &  &  &  &  &  &  &  &  &  & \cdots &  0\\
 &  &  &  &  &  &  &  &  &  &  &  & 1 & 1 & 1 & 1 \\
 &  &  &  &  &  &  &  &  &  &  &  &  &  &  & 0 \\
 \vdots&  &  &  &  &  &  &  &  &  &  &  &  &  &  & \vdots \\
 0&\cdots  &  &  &  &  &  &  &  &  &  &  &  &  & \cdots &0
\end{array}\right)
$$
from elementary row transformation, I found $\mathcal{C}_{16}:=\Lambda(\mathbb{F}_{2})$ is the vector space of degree $16$ and dimension $11$ over Finite Field of size 2 with the following generators. And the number of the elements is $2^{11} = 2048$ as the vector space is 11 dimensional over $\mathbb{F}_{2}$
$$
\mathcal{C}_{16}=\text{Span}_{\mathbb{F}_2}
\left(\begin{array}{llllllllllllllll}
1 &  &  &  &  &  &  & 1 &  &  &  & 1 &  & 1 & 1 & 1 \\
0 & 1 &  &  &  &  &  & 1 &  &  &  & 1 &  & 1 &  &  \\
\vdots &  & 1 &  &  &  &  & 1 &  &  &  & 1 &  &  & 1 &  \\
 &  &  & 1 &  &  &  & 1 &  &  &  & 1 &  &  &  & 1 \\
 &  &  &  & 1 &  &  & 1 &  &  &  &  &  & 1 & 1 &  \\
 &  &  &  &  & 1 &  & 1 &  &  &  &  &  & 1 &  & 1 \\
 &  &  &  &  &  & 1 & 1 &  &  &  &  &  &  & 1 & 1 \\
 &  &  &  &  &  &  &  & 1 &  &  & 1 &  & 1 & 1 &  \\
 &  &  &  &  &  &  &  &  & 1 &  & 1 &  & 1 &  & 1 \\
 \vdots&  &  &  &  &  &  &  &  &  & 1 & 1 &  &  & 1 & 1 \\
 0& \cdots &  &  &  &  &  &  &  &  & \cdots & 0 & 1 & 1 & 1 & 1
\end{array}\right)
$$
(Q.E.D)\\\\

Since all the number of the components, which is 1, are even of the generator of $\mathcal{C}_{16}$, the number of the components of the elements of it must be also even. And also by calculating all the combination of $\mathcal{C}_{16}$, I found the follwing result.
\begin{itembox}{Fact.1}
The elements of the linear space $\mathcal{C}_{16}$ has at least $4$ $1$'s as components ($140$ vectors). And the list of the  distribution is following.[See Appendix 5.3 for details.]
  $$
  \begin{array}{|c|c|c|c|c|c|c|}
  \hline 0 & 4 & 6 & 8 & 10 & 12 & 16 \\
  \hline 1 & 140 & 448 & 870 & 448 & 140 & 1 \\
  \hline
  \end{array}
  $$
\end{itembox}
Since $\pi$ is the homomorphism $\pi: \Lambda \rightarrow \mathcal{C}_{16}$. it can be also said
\begin{itembox}{Proposition.7}
For the type of the vectors $(2^n,-2^m,0^l)$ in $\Lambda$,
$$
n+m \geqslant 4
$$
where $\left(a^{n}, b^{m}\right)$ means a vector of type containing $n$ a's and $m$ b's.
\end{itembox}
(Proof)\\
If there exists the type $(2^n,-2^m,0^l)$ of the vectors such that $0<n+m<4$ in $\Lambda$, then $\mathcal{C}_{16}$ must have the vectors type $(1^k,0)$ such that $k<4$, as $\pi$ is the homomorphism. But actually not.Therefore, no such type of vector exists in $\Lambda$.(Q.E.D)\\\\
From this proposition, $\Lambda_{\mathbb{Q}^2_8}$ does not have the vectors, whose norm are $1/2$, such that
$$
\frac{1}{4}(2,2,0,0,...,0)
$$
By counting up the vectors that are the shortest of $\Lambda_{\mathbb{Q}^2_8}$, I find the following facts.
\begin{itembox}{Fact.2}
The lattice $\Lambda_{\mathbb{Q}^2_8}$ also has $4320$ shortest vectors.
\end{itembox}
The details of the above are as follows.
\begin{itembox}{Fact.3}
The vectors whose norm is 1 in $\Lambda_{\mathbb{Q}^2_8}$ are classified into the following types.
  $$
  \begin{array}{|l|l|}
  \hline\text{Types}&\text{Number of the vectors.}\\
  \hline\frac{1}{4}(-1^a,1^b)& 2048 \quad(a+b = 16) \\
  \frac{1}{4}\left(2^{4}, 0^{12}\right)& 140 \\
  \frac{1}{4}\left(-2^{4}, 0^{12}\right)& 140 \\
 \frac{1}{4}\left(4^{1}, 0^{15}\right)& 16 \\
  \frac{1}{4}\left(-4^{1}, 0^{15}\right)& 16 \\
  \frac{1}{4}\left(-2^{2}, 2^{2},0^{12}\right)& 840 \\
  \frac{1}{4}\left(-2^{1} \cdot 2^{3},0^{12}\right)&560 \\
  \frac{1}{4}\left(-2^{3}, 2^{1},0^{12}\right)& 560\\
  \hline
  \end{array}
  $$
where $(a^n,b^m)$ means a vector of type containing $n$ a's and $m$ b's.\\\\
Also vectors at the second distance are classified into the following types.
$$
\begin{array}{|l|l|}
\hline\frac{1}{4}\left(2^{6}, 0^{10}\right) & 448 \\
\frac{1}{4}\left(-2^{6}, 0^{10}\right) & 448 \\
\frac{1}{4}\left(-2^{1}, 2^{5}, 0^{10}\right) & 2688 \\
\frac{1}{4}\left(-2^{2}, 2^{4}, 0^{10}\right) & 6720 \\
\frac{1}{4}\left(-2^{3}, 2^{3}, 0^{10}\right) & 8960 \\
\frac{1}{4}\left(-2^{4}, 2^{2}, 0^{10}\right) & 6720 \\
\frac{1}{4}\left(-2^{5}, 2^{1}, 0^{10}\right) & 2688 \\
\frac{1}{4}\left(-3^{1}, 1^{a},-1^{b}\right) & 16384 \quad (a+b=15)\\
\frac{1}{4}\left(3^{1}, 1^{a},-1^{b}\right) & 16384 \quad (a+b=15)\\
\hline
\end{array}
$$
\end{itembox}
On the Other hand, the shortest vectors in $\Lambda_{\text{Barns-Wall}}$ are classified into follwoing types.(Total $4320$ vectors)
$$
\begin{array}{|l|l|}
\hline\text{Types} & \text{Number of the vectors}\\
\hline\frac{1}{4}\left(\pm 1^{8} , 0^{8}\right) & 30 \\
\frac{1}{4}\left(-1^{2} , 1^{6} , 0^{8}\right) & 30 \\
\frac{1}{4}\left(-1^{4} ,1^{4} , 0^{8}\right) & 30 \cdot\left(\begin{array}{l}
8 \\
2
\end{array}\right)=840 \\
\frac{1}{4}\left(-1^{4} , 1^{4} ,0^{8}\right) & 30 \cdot\left(\begin{array}{l}
8 \\
4
\end{array}\right)=2100 \\
\frac{1}{4}\left(-1^{6},1^{4} , 0^{8}\right) & 30 \cdot\left(\begin{array}{l}
8 \\
6
\end{array}\right)=840 \\
\frac{1}{4}\left(\pm 2^{2} , 0^{14}\right) & 240 \\
\frac{1}{4}\left(-2^{1} , 2^{1} ,0^{14}\right) & 2 \cdot 120=240 \\
\hline
\end{array}
$$

 Some vectors in the lattice $\Lambda_{\mathbb{Q}^2_8}$ with the elements in $G \subset \mathbb{Z}^{\prime}\left[\mathbb{Q}^2_8\right]$ are shown below.
$$
\begin{aligned}
&\frac{(1 + i + j + k)}{2} &\mapsto& (1/2, 1/2, 1/2, 1/2, 0, 0, 0, 0, 0, 0, 0, 0, 0, 0, 0, 0)\\
&\frac{(1 + i + j - k)}{2} &\mapsto &(1/2, 1/2, 1/2, -1/2, 0, 0, 0, 0, 0, 0, 0, 0, 0, 0, 0, 0)\\
&e_1=\frac{(1+i q+j r-k p)}{2} &\mapsto & (1/2, 0, 0, 0, 0, 0, 0, -1/2, 0, 1/2, 0, 0, 0, 0, 1/2, 0)\\
&e_{2}=\frac{(1-i p-j q-k r)}{2} &\mapsto & (1/2, 0, 0, 0, 0, -1/2, 0, 0, 0, 0, -1/2, 0, 0, 0, 0, -1/2)\\
&e_{3}=\frac{(1-i p-j r+k q)}{2}  &\mapsto & (1/2, 0, 0, 0, 0, -1/2, 0, 0, 0, 0, 0, 1/2, 0, 0, -1/2, 0)\\
&e_{4}=\frac{(1-i q-j p+k r)}{2} &\mapsto & (1/2, 0, 0, 1/2, 0, 0, 0, 0, -1/2, 0, 0, -1/2, 0, 0, 0, 0)\\
&\frac{(1 + i + j + k + p + ip + jp + \cdots + kr)}{4} &\mapsto & \frac{1}{4}(1, 1, 1, 1, 1, 1, 1, 1, 1, 1, 1, 1, 1, 1, 1, 1)
\end{aligned}
$$
I empirically predicted the following.
\begin{itembox}{Conjecture.1}
For the image $\text{Im}(\varphi^{-1})$ of the inverse mapping $\varphi^{-1}: G \rightarrow \{\mathbb{Z}\oplus(\mathbb{Z} + 1/2)\oplus (\mathbb{Z} + 1/4)\}^{16}$,
$$
\text{Im}(\varphi^{-1}) \subset \Lambda_{\mathbb{Q}^2_8}
$$
\end{itembox}

\section{Octanion and Weyl Group$\ W\left(E_{8}\right)$}
Octonion $\mathbb{O}=\left\{u_{0} \zeta_{0}+u_{1} \zeta_{1}+\cdots+u_{7} \zeta_{7} \mid u_{i} \in \mathbb{R}\right\}$is defined as $\zeta_i$ satisfying the following operations.
$$
\begin{array}{|c||c|c|c|c|c|c|c|c|}
\hline \times & \zeta_{0} & \zeta_{1} & \zeta_{2} & \zeta_{3} & \zeta_{4} & \zeta_{5} & \zeta_{6} & \zeta_{7} \\
\hline \hline\zeta_{0} & \zeta_{0} & \zeta_{1} & \zeta_{2} & \zeta_{3} & \zeta_{4} & \zeta_{5} & \zeta_{6} & \zeta_{7} \\
\hline \zeta_{1} & \zeta_{1} & -\zeta_{0} & \zeta_{3} & -\zeta_{2} & \zeta_{5} & -\zeta_{4} & -\zeta_{7} & \zeta_{6} \\
\hline \zeta_{2} & \zeta_{2} & -\zeta_{3} & -\zeta_{0} & \zeta_{1} & \zeta_{6} & \zeta_{7} & -\zeta_{4} & -\zeta_{5} \\
\hline \zeta_{3} & \zeta_{3} & \zeta_{2} & -\zeta_{1} & -\zeta_{0} & \zeta_{7} & -\zeta_{6} & \zeta_{5} & -\zeta_{4} \\
\hline \zeta_{4} & \zeta_{4} & -\zeta_{5} & -\zeta_{6} & -\zeta_{7} & -\zeta_{0} & \zeta_{1} & \zeta_{2} & \zeta_{3} \\
\hline \zeta_{5} & \zeta_{5} & \zeta_{4} & -\zeta_{7} & \zeta_{6} & -\zeta_{1} & -\zeta_{0} & -\zeta_{3} & \zeta_{2} \\
\hline \zeta_{6} & \zeta_{6} & \zeta_{7} & \zeta_{4} & -\zeta_{5} & -\zeta_{2} & \zeta_{3} & -\zeta_{0} & -\zeta_{1} \\
\hline \zeta_{7} & \zeta_{7} & -\zeta_{6} & \zeta_{5} & \zeta_{4} & -\zeta_{3} & -\zeta_{2} & \zeta_{1} & -\zeta_{0} \\
\hline
\end{array}
$$
The above rules of operation follow the arrows in the figure below.
\begin{figure}[h]
 \centering
 \includegraphics[keepaspectratio, scale=0.15]
      {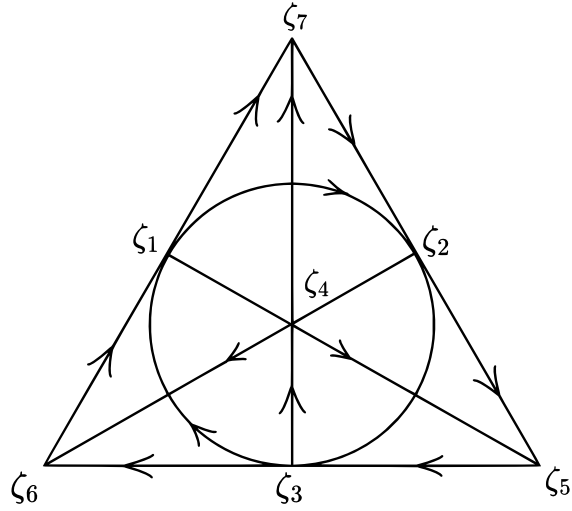}
 \caption{Figure representing the multiplication rules of Octantion.}
\label{}
\end{figure}

For example, an element$u \in \mathbb{O}$ can be written by the following linear combination.
$$
u=u_{0} \zeta_{0}+u_{1} \zeta_{1}+u_{2} \zeta_{2}+u_{3} \zeta_{3}+u_{4} \zeta_{4}+u_{5} \zeta_{5}+u_{6}\zeta_6 +u_{7} \zeta_{7}\\
$$
For the sake of simplicity, I will simply write the following.
$$
:=(u_0,u_1,u_2,u_3,u_4,u_5,u_6,u_7)
$$
then multiplying by $\zeta_1$ from the left will be
$$
\begin{aligned}
\zeta_{1} u &=-u_{1} \zeta_{}+u_{0} \zeta_{1}-u_{3} \zeta_{2}+u_{2} \zeta_{3}-u_{5} \zeta_{4}+u_{4} \zeta_{5}+u_{7} \zeta_{6}-u_{6} \zeta_{7} \\
&=(-u_1,u_0,-u_{3},u_{2},-u_{5},u_{4},u_{7},-u_{6})\\
&=\left(\begin{array}{cccc|cccc}
 & -1 &  &  &  &  &  &  \\
1 &  &  &  &  &  &  &  \\
 &  &  & -1 &  &  &  &  \\
 &  & 1 &  &  &  &  &  \\
\hline &  &  &  &  & -1 &  &  \\
 &  &  &  & 1 &  &  &  \\
 &  &  &  &  &  &  & 1 \\
 &  &  &  &  &  & -1 &
\end{array}\right)\left(\begin{array}{l}
u_{0} \\
u_{1} \\
u_{2} \\
u_{3} \\
u_{4} \\
u_{5} \\
u_{6} \\
u_{7}
\end{array}\right)
\end{aligned}
$$
Thus the above operation is the $8$-dimensional orthogonal transformation  $\zeta_{1} : \mathbb{R}^8 \rightarrow \mathbb{R}^8$  for $u = (u_0,u_1,...,u_7)$. From now on, I will write this $8 \times 8$ matrix simply as $\boldsymbol{\zeta}_{1}$.This can be written using the matrices described in the previous section
$$
\boldsymbol{\zeta}_{1}=\left(\begin{array}{c|c}
-\boldsymbol{i} &  \\
\hline& -\boldsymbol{p}
\end{array}\right)
$$
Others can also be written similarly as follows
$$
\begin{array}{l}
\boldsymbol{\zeta}_{2}=\left(\begin{array}{c|c}
-\boldsymbol{j} &  \\
\hline  & -\boldsymbol{q}
\end{array}\right) \\\\
\boldsymbol{\zeta}_{3}=\left(\begin{array}{c|c}
\boldsymbol{k} &  \\
\hline  & \boldsymbol{r}
\end{array}\right)\\\\
\boldsymbol{\zeta}_{4}=\left(\begin{array}{c|c}
 & \frac{1}{2}(\boldsymbol{i} \boldsymbol{p}+\boldsymbol{j} \boldsymbol{q}-\boldsymbol{k} \boldsymbol{r}+E) \\
\hline -\frac{1}{2}(\boldsymbol{i} \boldsymbol{p}+\boldsymbol{j} \boldsymbol{q}-\boldsymbol{k} \boldsymbol{r}+E) &
\end{array}\right)\\\\
\boldsymbol{\zeta}_{5}=\left(\begin{array}{c|c}
 & \frac{1}{2}(\boldsymbol{i}-\boldsymbol{jr}- \boldsymbol{kq}-\boldsymbol{p}) \\
\hline \frac{1}{2}(\boldsymbol{i}+\boldsymbol{jr} +\boldsymbol{kq}-\boldsymbol{p}) &
\end{array}\right)\\\\
\boldsymbol{\zeta}_{6}=\left(\begin{array}{c|c}
 & \frac{1}{2}(\boldsymbol{ir}+\boldsymbol{j}+\boldsymbol{kp}-\boldsymbol{q}) \\
\hline \frac{1}{2}(-\boldsymbol{ir}+\boldsymbol{j}-\boldsymbol{kp}-\boldsymbol{q}) &
\end{array}\right)\\\\
\boldsymbol{\zeta}_7 = \left(\begin{array}{c|c}
 & \frac{1}{2}(-\boldsymbol{iq}+\boldsymbol{jp}-\boldsymbol{k}-\boldsymbol{q}) \\
\hline \frac{1}{2}(\boldsymbol{iq}-\boldsymbol{jp}-\boldsymbol{k}-\boldsymbol{q}) &
\end{array}\right)
\end{array}
$$
Now I define
$$
\boldsymbol{\eta}:=\frac{\boldsymbol{\zeta}_{1}+\boldsymbol{\zeta}_{2}+\boldsymbol{\zeta}_{3}+\boldsymbol{\zeta}_{4}}{2}=\left(\begin{array}{c|c}
\frac{1}{2}(-\boldsymbol{i}-\boldsymbol{j}+\boldsymbol{k}) & \frac{1}{4}(\boldsymbol{ip}+\boldsymbol{jq}-\boldsymbol{kr}+E) \\
\hline-\frac{1}{4}(\boldsymbol{ip}+\boldsymbol{jq}-\boldsymbol{kr}+E)  & \frac{1}{2}(-\boldsymbol{p}-\boldsymbol{q}+\boldsymbol{r})
\end{array}\right)
$$
And
$$
\begin{array}{l}
\boldsymbol{\xi}_{1}:=\left(\begin{array}{c|c}
-\boldsymbol{p} &  \\
\hline  & \boldsymbol{p}
\end{array}\right) \\
\boldsymbol{\xi}_{2}:=\left(\begin{array}{c|c}
-\boldsymbol{q} &  \\
\hline  & \boldsymbol{q}
\end{array}\right) \\
\boldsymbol{\xi}_{3}:=\left(\begin{array}{c|c}
-\boldsymbol{r} &  \\
\hline  & \boldsymbol{r}
\end{array}\right)
\end{array}
$$
then I observed the following.
\begin{itembox}{Observation}
  The group $G^{ \prime}$ defined below is the orthogonal group of order $348364800 = 2^{13} \cdot 3^{5} \cdot 5^{2} \cdot 7$
  $$
  \begin{aligned}
  G^{ \prime} &=\left\langle\boldsymbol{\zeta}_{1}, \boldsymbol{\zeta}_{2}, \boldsymbol{\zeta}_{3}, \boldsymbol{\eta}, \boldsymbol{\xi}_{1}, \boldsymbol{\xi}_{2}\right\rangle
  \end{aligned}
  $$
  and satisfies
  $$
 \simeq O_{8}^{+}(2) =W\left(E_{8}\right)
  $$
\end{itembox}

\section{Appendex}
\subparagraph{matrixies}
$$
\begin{array}{l}
{e}_{1}=\frac{1}{2}(1, i, j, k)\left(\begin{array}{cccc}
1 & 0 & 0 & 0 \\
0 & 0 & 1 & 0 \\
0 & 0 & 0 & 1 \\
0 & -1 & 0 & 0
\end{array}\right)\left(\begin{array}{l}
1 \\
p \\
q \\
r
\end{array}\right) \stackrel{\rho}{\longmapsto} \frac{1}{2}\left(\begin{array}{cccc}
1 & 1 & 1 & 1 \\
1 & 1 & -1 & -1 \\
1 & -1 & 1 & -1 \\
1 & -1 & -1 & 1
\end{array}\right)\\
e_{2}=\frac{1}{2}(1, i, j, k)\left(\begin{array}{cccc}
1 & 0 & 0 & 0 \\
0 & -1 & 0 & 0 \\
0 & 0 & -1 & 0 \\
0 & 0 & 0 & -1
\end{array}\right)\left(\begin{array}{l}
1 \\
p \\
q \\
r
\end{array}\right)\stackrel{\rho}{\longmapsto}\left(\begin{array}{rrrr}
1 & 0 & 0 & 0 \\
0 & 1 & 0 & 0 \\
0 & 0 & 1 & 0 \\
0 & 0 & 0 & -1
\end{array}\right)\\
e_{3}=\frac{1}{2}(1, i, j, k)\left(\begin{array}{cccc}
1 & 0 & 0 & 0 \\
0 & -1 & 0 & 0 \\
0 & 0 & 0 & -1 \\
0 & 0 & 1 & 0
\end{array}\right)\left(\begin{array}{l}
1 \\
p \\
q \\
r
\end{array}\right)\stackrel{\rho}{\longmapsto}\left(\begin{array}{llll}
1 & 0 & 0 & 0 \\
0 & 1 & 0 & 0 \\
0 & 0 & 0 & 1 \\
0 & 0 & 1 & 0
\end{array}\right)\\
e_{4}=\frac{1}{2}(1, i, j, k)\left(\begin{array}{cccc}
1 & 0 & 0 & 0 \\
0 & 0 & -1 & 0 \\
0 & -1 & 0 & 0 \\
0 & 0 & 0 & 1
\end{array}\right)\left(\begin{array}{l}
1 \\
p \\
q \\
r
\end{array}\right)\stackrel{\rho}{\longmapsto} \left(\begin{array}{llll}
1 & 0 & 0 & 0 \\
0 & 0 & 1 & 0 \\
0 & 1 & 0 & 0 \\
0 & 0 & 0 & 1
\end{array}\right)
\end{array}
$$

$$
\boldsymbol{\zeta}_1=\left(\begin{array}{rrrrrrrr}
0 & -1 & 0 & 0 & 0 & 0 & 0 & 0 \\
1 & 0 & 0 & 0 & 0 & 0 & 0 & 0 \\
0 & 0 & 0 & -1 & 0 & 0 & 0 & 0 \\
0 & 0 & 1 & 0 & 0 & 0 & 0 & 0 \\
0 & 0 & 0 & 0 & 0 & -1 & 0 & 0 \\
0 & 0 & 0 & 0 & 1 & 0 & 0 & 0 \\
0 & 0 & 0 & 0 & 0 & 0 & 0 & 1 \\
0 & 0 & 0 & 0 & 0 & 0 & -1 & 0
\end{array}\right)
$$

$$
\boldsymbol{\zeta}_2=\left(\begin{array}{rrrrrrrr}
0 & 0 & -1 & 0 & 0 & 0 & 0 & 0 \\
0 & 0 & 0 & 1 & 0 & 0 & 0 & 0 \\
1 & 0 & 0 & 0 & 0 & 0 & 0 & 0 \\
0 & -1 & 0 & 0 & 0 & 0 & 0 & 0 \\
0 & 0 & 0 & 0 & 0 & 0 & -1 & 0 \\
0 & 0 & 0 & 0 & 0 & 0 & 0 & -1 \\
0 & 0 & 0 & 0 & 1 & 0 & 0 & 0 \\
0 & 0 & 0 & 0 & 0 & 1 & 0 & 0
\end{array}\right)
$$

$$
\boldsymbol{\zeta}_3=\left(\begin{array}{rrrrrrrr}
0 & 0 & 0 & -1 & 0 & 0 & 0 & 0 \\
0 & 0 & -1 & 0 & 0 & 0 & 0 & 0 \\
0 & 1 & 0 & 0 & 0 & 0 & 0 & 0 \\
1 & 0 & 0 & 0 & 0 & 0 & 0 & 0 \\
0 & 0 & 0 & 0 & 0 & 0 & 0 & -1 \\
0 & 0 & 0 & 0 & 0 & 0 & 1 & 0 \\
0 & 0 & 0 & 0 & 0 & -1 & 0 & 0 \\
0 & 0 & 0 & 0 & 1 & 0 & 0 & 0
\end{array}\right)
$$

$$
\boldsymbol{\zeta}_4=\left(\begin{array}{rrrrrrrr}
0 & 0 & 0 & 0 & -1 & 0 & 0 & 0 \\
0 & 0 & 0 & 0 & 0 & 1 & 0 & 0 \\
0 & 0 & 0 & 0 & 0 & 0 & 1 & 0 \\
0 & 0 & 0 & 0 & 0 & 0 & 0 & 1 \\
1 & 0 & 0 & 0 & 0 & 0 & 0 & 0 \\
0 & -1 & 0 & 0 & 0 & 0 & 0 & 0 \\
0 & 0 & -1 & 0 & 0 & 0 & 0 & 0 \\
0 & 0 & 0 & -1 & 0 & 0 & 0 & 0
\end{array}\right)
$$

$$
\boldsymbol{\zeta}_5=\left(\begin{array}{rrrrrrrr}
0 & 0 & 0 & 0 & 0 & -1 & 0 & 0 \\
0 & 0 & 0 & 0 & -1 & 0 & 0 & 0 \\
0 & 0 & 0 & 0 & 0 & 0 & 0 & 1 \\
0 & 0 & 0 & 0 & 0 & 0 & -1 & 0 \\
0 & 1 & 0 & 0 & 0 & 0 & 0 & 0 \\
1 & 0 & 0 & 0 & 0 & 0 & 0 & 0 \\
0 & 0 & 0 & 1 & 0 & 0 & 0 & 0 \\
0 & 0 & -1 & 0 & 0 & 0 & 0 & 0
\end{array}\right)
$$

$$
\boldsymbol{\zeta}_6=\left(\begin{array}{rrrrrrrr}
0 & 0 & 0 & 0 & 0 & 0 & -1 & 0 \\
0 & 0 & 0 & 0 & 0 & 0 & 0 & -1 \\
0 & 0 & 0 & 0 & -1 & 0 & 0 & 0 \\
0 & 0 & 0 & 0 & 0 & 1 & 0 & 0 \\
0 & 0 & 1 & 0 & 0 & 0 & 0 & 0 \\
0 & 0 & 0 & -1 & 0 & 0 & 0 & 0 \\
1 & 0 & 0 & 0 & 0 & 0 & 0 & 0 \\
0 & 1 & 0 & 0 & 0 & 0 & 0 & 0
\end{array}\right)
$$
$$
\boldsymbol{\zeta}_7=\left(\begin{array}{rrrrrrrr}
0 & 0 & 0 & 0 & 0 & 0 & 0 & -1 \\
0 & 0 & 0 & 0 & 0 & 0 & 1 & 0 \\
0 & 0 & 0 & 0 & 0 & -1 & 0 & 0 \\
0 & 0 & 0 & 0 & -1 & 0 & 0 & 0 \\
0 & 0 & 0 & 1 & 0 & 0 & 0 & 0 \\
0 & 0 & 1 & 0 & 0 & 0 & 0 & 0 \\
0 & -1 & 0 & 0 & 0 & 0 & 0 & 0 \\
1 & 0 & 0 & 0 & 0 & 0 & 0 & 0
\end{array}\right)
$$

$$
\boldsymbol{\eta}=\frac{1}{2}\left(\begin{array}{rrrrrrrr}
0 & 1 & 1 & 1 & 1 & 0 & 0 & 0 \\
-1 & 0 & 1 & -1 & 0 & -1 & 0 & 0 \\
-1 & -1 & 0 & 1 & 0 & 0 & -1 & 0 \\
-1 & 1 & -1 & 0 & 0 & 0 & 0 & -1 \\
-1 & 0 & 0 & 0 & 0 & 1 & 1 & 1 \\
0 & 1 & 0 & 0 & -1 & 0 & -1 & 1 \\
0 & 0 & 1 & 0 & -1 & 1 & 0 & -1 \\
0 & 0 & 0 & 1 & -1 & -1 & 1 & 0
\end{array}\right)
$$

$$
\boldsymbol{\xi}_1=\left(\begin{array}{rrrrrrrr}
0 & -1 & 0 & 0 & 0 & 0 & 0 & 0 \\
1 & 0 & 0 & 0 & 0 & 0 & 0 & 0 \\
0 & 0 & 0 & 1 & 0 & 0 & 0 & 0 \\
0 & 0 & -1 & 0 & 0 & 0 & 0 & 0 \\
0 & 0 & 0 & 0 & 0 & 1 & 0 & 0 \\
0 & 0 & 0 & 0 & -1 & 0 & 0 & 0 \\
0 & 0 & 0 & 0 & 0 & 0 & 0 & -1 \\
0 & 0 & 0 & 0 & 0 & 0 & 1 & 0
\end{array}\right)
$$

$$
\boldsymbol{\xi}_2=\left(\begin{array}{rrrrrrrr}
0 & 0 & -1 & 0 & 0 & 0 & 0 & 0 \\
0 & 0 & 0 & -1 & 0 & 0 & 0 & 0 \\
1 & 0 & 0 & 0 & 0 & 0 & 0 & 0 \\
0 & 1 & 0 & 0 & 0 & 0 & 0 & 0 \\
0 & 0 & 0 & 0 & 0 & 0 & 1 & 0 \\
0 & 0 & 0 & 0 & 0 & 0 & 0 & 1 \\
0 & 0 & 0 & 0 & -1 & 0 & 0 & 0 \\
0 & 0 & 0 & 0 & 0 & -1 & 0 & 0
\end{array}\right)
$$

$$
\boldsymbol{\xi}_3=\left(\begin{array}{rrrrrrrr}
0 & 0 & 0 & -1 & 0 & 0 & 0 & 0 \\
0 & 0 & 1 & 0 & 0 & 0 & 0 & 0 \\
0 & -1 & 0 & 0 & 0 & 0 & 0 & 0 \\
1 & 0 & 0 & 0 & 0 & 0 & 0 & 0 \\
0 & 0 & 0 & 0 & 0 & 0 & 0 & 1 \\
0 & 0 & 0 & 0 & 0 & 0 & -1 & 0 \\
0 & 0 & 0 & 0 & 0 & 1 & 0 & 0 \\
0 & 0 & 0 & 0 & -1 & 0 & 0 & 0
\end{array}\right)
$$
\subsection{The permutation representation of $\mathbb{Q}^2_8$}
By the way, we can find that $i,j,k,p,q,r$ as a subgroup of the group $C \subset Sym(16)$ order $192$ generated by the following permutations of the rotations in three orthogonal directions.
$$
\begin{aligned}
C=&<(1,2,10,9)(3,11,12,4)(5,6,14,13)(7,15,16,8), \\
&(1,4,12,9)(2,3,11,10)(5,8,16,13)(6,7,15,14), \\
&(1,5,13,9)(2,6,14,10)(3,7,15,11)(4,8,16,12)>
\end{aligned}
$$
\begin{figure}[h]
 \centering
 \includegraphics[keepaspectratio, scale=0.3]{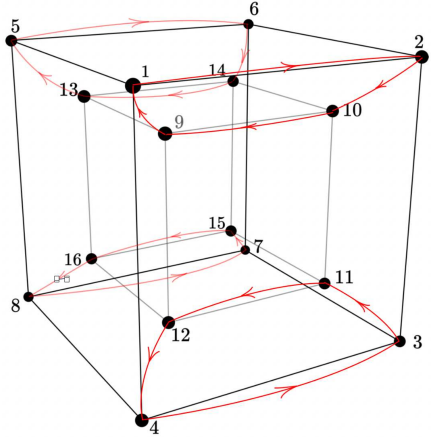}
\caption{The permutation of the generator $g_1$ of the group $C$}
\label{}
\end{figure}

The figure shows  one of the generator which is
$$
g_1 = (1,2,10,9)(3,11,12,4)(5,6,14,13)(7,15,16,8) \in C
$$
and $C$ include subset $C^\prime$ order $32$ generated by
$$
\begin{aligned}
&(1,3,15,13)(2,11,16,5)(4,7,14,9)(6,10,12,8)\\
&(1,12,15,6)(2,9,16,7)(3,10,13,8)(4,11,14,5)\\
&(1,8,15,10)(2,5,16,11)(3,6,13,12)(4,7,14,9)\\
&(1,13,15,3),(2,9,16,7),(4,5,14,11),(6,10,12,8)
\end{aligned}
$$
and there is isomorphism $\tau:\mathbb{Q}^2_8 \cong <i,j,p,q> \rightarrow C^\prime$ that is
$$
\begin{aligned}
&(i)\tau=(1,3,15,13)(2,11,16,5)(4,7,14,9)(6,10,12,8) \\
&(j)\tau=(1,12,15,6)(2,9,16,7)(3,10,13,8)(4,11,14,5) \\
&(k)\tau=(1,10,15,8)(2,14,16,4)(3,6,13,12)(5,9,11,7)\\
&(p)\tau=(1,8,15,10)(2,5,16,11)(3,6,13,12)(4,7,14,9) \\
&(q)\tau=(1,13,15,3),(2,9,16,7),(4,5,14,11),(6,10,12,8)\\
&(r)\tau=(1,6,15,12)(2,14,16,4)(3,10,13,8)(5,7,11,9)\\
&(-1)\tau=(1,15)(2,16)(3,13)(4,14)(5,11)(6,12)(7,9)(8,10)
\end{aligned}
$$
As $-1$ is commutative with all the elements, actually $(-1)\tau=(1,15)(2,16)(3,13)(4,14)(5,11)(6,12)(7,9)(8,10)$ is the center of the group $C$. Graphically
\begin{figure}[h]
 \centering
 \includegraphics[keepaspectratio, scale=0.4]
      {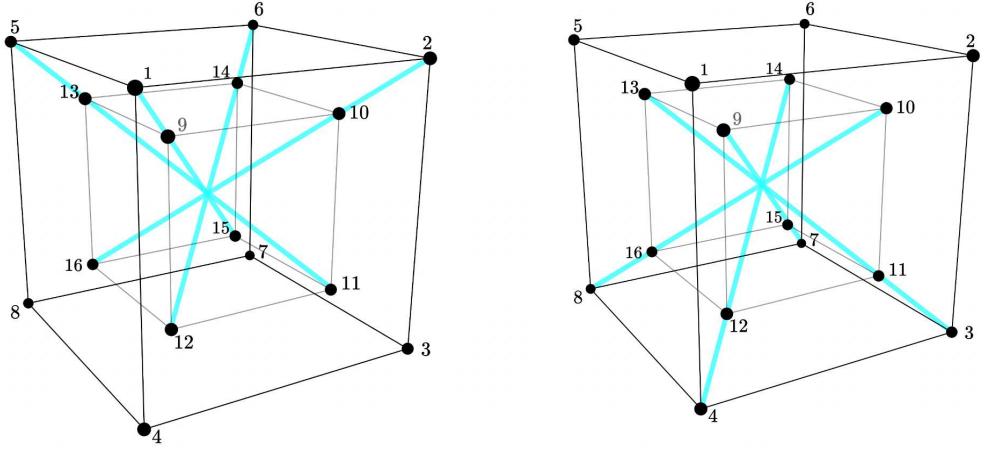}
 \caption{The center of $C$}
\label{}
\end{figure}

\subsection{A constraction of Quotanion}
$$
\mathcal{I} = (\sqrt{-1}) \sigma=\left(\begin{array}{cc}
 & 1 \\
-1 &
\end{array}\right)
$$
and let $ \mathcal{J}$ to be
$$
\mathcal{I} \mathcal{J} = - \mathcal{J} \mathcal{I}
$$
then $ \mathcal{J}$ is determined as
$$
\mathcal{J}=\left(\begin{array}{ll}
1 &0 \\
0& -1
\end{array}\right)
$$
$\mathcal{I}$ and $\mathcal{J}$ satisfy
$$
\mathcal{I}^2 = -1,\quad \mathcal{I}^4 = \mathcal{J}^2 = 1,\quad \mathcal{I} \mathcal{J}=-\mathcal{J} \mathcal{I}
$$

using above
$$
\begin{array}{l}
\boldsymbol{i}=\left(\begin{array}{ll}
\mathcal{J} & \\
& \mathcal{I}
\end{array}\right) \quad \boldsymbol{k}=\left(\begin{array}{ll}
& -1 \\
1 &
\end{array}\right)\\
\boldsymbol{j}=\left(\begin{array}{cc}
& \mathcal{I}\\
-\mathcal{I} &
\end{array}\right)
\end{array}
$$
Also the multiplicative group generated by $\mathcal{I}, \mathcal{J}$ order $8$ is isomorphic to $W(B_2) \cong W(C_2)$

\subsection{The code $\mathcal{C}_{16}$ in $(\mathbb{F}_2)^{16}$}
The generator matrix of $\mathcal{C}_{16}$ is
$$
\left(\begin{array}{c}
u_{1} \\
u_{2} \\
u_{3} \\
\vdots \\
 \\
 \\
 \\
 \\
 \vdots\\
u_{11}
\end{array}\right)=
\left(\begin{array}{llllllllllllllll}
1 & 0 & 0 & 0 & 0 & 0 & 0 & 1 & 0 & 0 & 0 & 1 & 0 & 1 & 1 & 1 \\
0 & 1 & 0 & 0 & 0 & 0 & 0 & 1 & 0 & 0 & 0 & 1 & 0 & 1 & 0 & 0 \\
0 & 0 & 1 & 0 & 0 & 0 & 0 & 1 & 0 & 0 & 0 & 1 & 0 & 0 & 1 & 0 \\
0 & 0 & 0 & 1 & 0 & 0 & 0 & 1 & 0 & 0 & 0 & 1 & 0 & 0 & 0 & 1 \\
0 & 0 & 0 & 0 & 1 & 0 & 0 & 1 & 0 & 0 & 0 & 0 & 0 & 1 & 1 & 0 \\
0 & 0 & 0 & 0 & 0 & 1 & 0 & 1 & 0 & 0 & 0 & 0 & 0 & 1 & 0 & 1 \\
0 & 0 & 0 & 0 & 0 & 0 & 1 & 1 & 0 & 0 & 0 & 0 & 0 & 0 & 1 & 1 \\
0 & 0 & 0 & 0 & 0 & 0 & 0 & 0 & 1 & 0 & 0 & 1 & 0 & 1 & 1 & 0 \\
0 & 0 & 0 & 0 & 0 & 0 & 0 & 0 & 0 & 1 & 0 & 1 & 0 & 1 & 0 & 1 \\
0 & 0 & 0 & 0 & 0 & 0 & 0 & 0 & 0 & 0 & 1 & 1 & 0 & 0 & 1 & 1 \\
0 & 0 & 0 & 0 & 0 & 0 & 0 & 0 & 0 & 0 & 0 & 0 & 1 & 1 & 1 & 1
\end{array}\right)
$$
For all vectors $x \in \mathcal{C}_{16}$ can be written as
$$
x = a_1u_1+a_2u_2 + \cdots +a_{11}u_{11}
$$
where $a_i \in \mathbb{F}_{2}$.The operation is then performed as follows.
$$
\begin{aligned}
u_{10} + u_{11} &= (0,0,\cdots,0,1,1,0,0,1,1) + (0,0,\cdots,0,0,0,1,1,1,1)\\
&=(0,0,\cdots,0,1,1,1,1,0,0)
\end{aligned}
$$
This generates the vector space of degree $16$ and dimension $11$ over Finite Field of size $2$.\\\\
The weight distribution of $\mathcal{C}_{16}$ is as follows.
$$
\begin{array}{|c|c|c|c|c|c|c|}
\hline 0 & 4 & 6 & 8 & 10 & 12 & 16 \\
\hline 1 & 140 & 448 & 870 & 448 & 140 & 1 \\
\hline
\end{array}
$$
and at least, $\mathcal{C}_{16}$ is invariant by the action of the group $G_{\mathcal{C}}$ order $322560=2^{10} \cdot 3^{2} \cdot 5 \cdot 7$ defined by the following
$$
G_\mathcal{C}:=\langle(0,1)(2,7,9,11,13,4,14,3,6,8,10,12,5,15),(0,15,8,11,14,1,6,5)(2,3,10,7,12,13,4,9)\rangle
$$
and I also found $G_\mathcal{C} \cong 2^{4}\cdot A_{8}$ which is called 'Miracle Octad Group' of Mathieu group $M_{24}$, the subgroup of $G_{\mathcal{C}}$ generated by following is isomorphic to $A_8$.
$$
\begin{aligned}
\tau_1=&(1,4,2,7)(5,6)(8,13,11,14)(12,15)\\
\tau_2=&(1,6)(2,14)(3,8)(4,15)(5,9)(10,13)
\end{aligned}
$$
 \\\\

For any elements of $\mathcal{C}_{16}$(which has four $1$s) share $2$ coordinates for each of the $32$ generators(which has four $1$s), and $1$ coordinates for each of the $64$ generators(which has four $1$s). Below are 140 vectors in $\mathcal{C}_{16}$ with four 1's in the coordinates.
$$
\begin{aligned}
&(0,1,0,0,0,0,0,1,0,0,0,1,0,1,0,0) \\
&(1,1,0,0,0,0,0,0,0,0,0,0,0,0,1,1) \\
&(0,0,1,0,0,0,0,1,0,0,0,1,0,0,1,0) \\
&(1,0,1,0,0,0,0,0,0,0,0,0,0,1,0,1) \\
&(0,1,1,0,0,0,0,0,0,0,0,0,0,1,1,0) \\
&(0,0,0,1,0,0,0,1,0,0,0,1,0,0,0,1) \\
&(1,0,0,1,0,0,0,0,0,0,0,0,0,1,1,0) \\
&(0,1,0,1,0,0,0,0,0,0,0,0,0,1,0,1) \\
&(0,0,1,1,0,0,0,0,0,0,0,0,0,0,1,1) \\
&(1,1,1,1,0,0,0,0,0,0,0,0,0,0,0,0)
\end{aligned}
$$
$$
\begin{aligned}
&(0,0,0,0,1,0,0,1,0,0,0,0,0,1,1,0) \\
&(1,0,0,0,1,0,0,0,0,0,0,1,0,0,0,1) \\
&(0,1,0,0,1,0,0,0,0,0,0,1,0,0,1,0) \\
&(0,0,1,0,1,0,0,0,0,0,0,1,0,1,0,0) \\
&(0,1,1,0,1,0,0,1,0,0,0,0,0,0,0,0) \\
&(1,0,0,1,1,0,0,1,0,0,0,0,0,0,0,0) \\
&(0,0,0,0,0,1,0,1,0,0,0,0,0,1,0,1) \\
&(1,0,0,0,0,1,0,0,0,0,0,1,0,0,1,0) \\
&(0,1,0,0,0,1,0,0,0,0,0,1,0,0,0,1) \\
&(1,0,1,0,0,1,0,1,0,0,0,0,0,0,0,0)
\end{aligned}
$$
$$
\begin{aligned}
&(0,0,0,1,0,1,0,0,0,0,0,1,0,1,0,0) \\
&(0,1,0,1,0,1,0,1,0,0,0,0,0,0,0,0) \\
&(0,0,0,0,1,1,0,0,0,0,0,0,0,0,1,1) \\
&(1,1,0,0,1,1,0,0,0,0,0,0,0,0,0,0) \\
&(0,0,1,1,1,1,0,0,0,0,0,0,0,0,0,0) \\
&(0,0,0,0,0,0,1,1,0,0,0,0,0,0,1,1) \\
&(1,0,0,0,0,0,1,0,0,0,0,1,0,1,0,0) \\
&(1,1,0,0,0,0,1,1,0,0,0,0,0,0,0,0) \\
&(0,0,1,0,0,0,1,0,0,0,0,1,0,0,0,1) \\
&(0,0,0,1,0,0,1,0,0,0,0,1,0,0,1,0)
\end{aligned}
$$
$$
\begin{aligned}
&(0,0,1,1,0,0,1,1,0,0,0,0,0,0,0,0) \\
&(0,0,0,0,1,0,1,0,0,0,0,0,0,1,0,1) \\
&(1,0,1,0,1,0,1,0,0,0,0,0,0,0,0,0) \\
&(0,1,0,1,1,0,1,0,0,0,0,0,0,0,0,0) \\
&(0,0,0,0,0,1,1,0,0,0,0,0,0,1,1,0) \\
&(0,1,1,0,0,1,1,0,0,0,0,0,0,0,0,0) \\
&(1,0,0,1,0,1,1,0,0,0,0,0,0,0,0,0) \\
&(0,0,0,0,1,1,1,1,0,0,0,0,0,0,0,0) \\
&(0,0,0,0,0,0,0,0,1,0,0,1,0,1,1,0) \\
&(1,0,0,0,0,0,0,1,1,0,0,0,0,0,0,1)
\end{aligned}
$$
$$
\begin{aligned}
&(0,1,0,0,0,0,0,1,1,0,0,0,0,0,1,0) \\
&(0,0,1,0,0,0,0,1,1,0,0,0,0,1,0,0) \\
&(0,1,1,0,0,0,0,0,1,0,0,1,0,0,0,0) \\
&(1,0,0,1,0,0,0,0,1,0,0,1,0,0,0,0) \\
&(0,0,0,0,1,0,0,1,1,0,0,1,0,0,0,0) \\
&(0,1,0,0,1,0,0,0,1,0,0,0,0,1,0,0) \\
&(0,0,1,0,1,0,0,0,1,0,0,0,0,0,1,0) \\
&(0,0,0,1,1,0,0,0,1,0,0,0,0,0,0,1) \\
&(1,0,0,0,0,1,0,0,1,0,0,0,0,1,0,0) \\
&(0,0,1,0,0,1,0,0,1,0,0,0,0,0,0,1)
\end{aligned}
$$
$$
\begin{aligned}
&(0,0,0,1,0,1,0,0,1,0,0,0,0,0,1,0) \\
&(1,0,0,0,0,0,1,0,1,0,0,0,0,0,1,0) \\
&(0,1,0,0,0,0,1,0,1,0,0,0,0,0,0,1) \\
&(0,0,0,1,0,0,1,0,1,0,0,0,0,1,0,0) \\
&(0,0,0,0,0,1,1,0,1,0,0,1,0,0,0,0) \\
&(0,0,0,0,0,0,0,0,0,1,0,1,0,1,0,1) \\
&(1,0,0,0,0,0,0,1,0,1,0,0,0,0,1,0) \\
&(0,1,0,0,0,0,0,1,0,1,0,0,0,0,0,1) \\
&(1,0,1,0,0,0,0,0,0,1,0,1,0,0,0,0) \\
&(0,0,0,1,0,0,0,1,0,1,0,0,0,1,0,0)
\end{aligned}
$$
$$
\begin{aligned}
&(0,1,0,1,0,0,0,0,0,1,0,1,0,0,0,0) \\
&(1,0,0,0,1,0,0,0,0,1,0,0,0,1,0,0) \\
&(0,0,1,0,1,0,0,0,0,1,0,0,0,0,0,1) \\
&(0,0,0,1,1,0,0,0,0,1,0,0,0,0,1,0) \\
&(0,0,0,0,0,1,0,1,0,1,0,1,0,0,0,0) \\
&(0,1,0,0,0,1,0,0,0,1,0,0,0,1,0,0) \\
&(0,0,1,0,0,1,0,0,0,1,0,0,0,0,1,0) \\
&(0,0,0,1,0,1,0,0,0,1,0,0,0,0,0,1) \\
&(1,0,0,0,0,0,1,0,0,1,0,0,0,0,0,1) \\
&(0,1,0,0,0,0,1,0,0,1,0,0,0,0,1,0)
\end{aligned}
$$
$$
\begin{aligned}
&(0,0,1,0,0,0,0,1,0,0,1,0,0,0,0,1) \\
&(0,0,0,1,0,0,0,1,0,0,1,0,0,0,1,0) \\
&(0,0,1,1,0,0,0,0,0,0,1,1,0,0,0,0) \\
&(1,0,0,0,1,0,0,0,0,0,1,0,0,0,1,0) \\
&(0,1,0,0,1,0,0,0,0,0,1,0,0,0,0,1) \\
&(0,0,0,1,1,0,0,0,0,0,1,0,0,1,0,0) \\
&(1,0,0,0,0,1,0,0,0,0,1,0,0,0,0,1) \\
&(0,1,0,0,0,1,0,0,0,0,1,0,0,0,1,0) \\
&(0,0,1,0,0,1,0,0,0,0,1,0,0,1,0,0) \\
&(0,0,0,0,1,1,0,0,0,0,1,1,0,0,0,0)
\end{aligned}
$$
$$
\begin{aligned}
&(0,0,0,0,0,0,1,1,0,0,1,1,0,0,0,0) \\
&(0,1,0,0,0,0,1,0,0,0,1,0,0,1,0,0) \\
&(0,0,1,0,0,0,1,0,0,0,1,0,0,0,1,0) \\
&(0,0,0,1,0,0,1,0,0,0,1,0,0,0,0,1) \\
&(0,0,0,0,0,0,0,0,1,0,1,0,0,1,0,1) \\
&(1,0,1,0,0,0,0,0,1,0,1,0,0,0,0,0) \\
&(0,1,0,1,0,0,0,0,1,0,1,0,0,0,0,0) \\
&(0,0,0,0,0,1,0,1,1,0,1,0,0,0,0,0) \\
&(0,0,0,0,1,0,1,0,1,0,1,0,0,0,0,0) \\
&(0,0,0,0,0,0,0,0,0,1,1,0,0,1,1,0)
\end{aligned}
$$
$$
\begin{aligned}
&(0,1,1,0,0,0,0,0,0,1,1,0,0,0,0,0) \\
&(1,0,0,1,0,0,0,0,0,1,1,0,0,0,0,0) \\
&(0,0,0,0,1,0,0,1,0,1,1,0,0,0,0,0) \\
&(0,0,0,0,0,1,1,0,0,1,1,0,0,0,0,0) \\
&(0,0,0,0,0,0,0,0,1,1,1,1,0,0,0,0) \\
&(0,0,0,0,0,0,0,0,0,0,0,0,1,1,1,1) \\
&(1,0,0,0,0,0,0,1,0,0,0,1,1,0,0,0) \\
&(1,1,0,0,0,0,0,0,0,0,0,0,1,1,0,0) \\
&(1,0,1,0,0,0,0,0,0,0,0,0,1,0,1,0) \\
&(0,1,1,0,0,0,0,0,0,0,0,0,1,0,0,1)
\end{aligned}
$$
$$
\begin{aligned}
&(1,0,0,1,0,0,0,0,0,0,0,0,1,0,0,1) \\
&(0,1,0,1,0,0,0,0,0,0,0,0,1,0,1,0) \\
&(0,0,1,1,0,0,0,0,0,0,0,0,1,1,0,0) \\
&(0,0,0,0,1,0,0,1,0,0,0,0,1,0,0,1) \\
&(0,0,0,1,1,0,0,0,0,0,0,1,1,0,0,0) \\
&(0,0,0,0,0,1,0,1,0,0,0,0,1,0,1,0) \\
&(0,0,1,0,0,1,0,0,0,0,0,1,1,0,0,0) \\
&(0,0,0,0,1,1,0,0,0,0,0,0,1,1,0,0) \\
&(0,0,0,0,0,0,1,1,0,0,0,0,1,1,0,0) \\
&(0,1,0,0,0,0,1,0,0,0,0,1,1,0,0,0)
\end{aligned}
$$
$$
\begin{aligned}
&(0,0,0,0,1,0,1,0,0,0,0,0,1,0,1,0) \\
&(0,0,0,0,0,1,1,0,0,0,0,0,1,0,0,1) \\
&(0,0,0,0,0,0,0,0,1,0,0,1,1,0,0,1) \\
&(0,0,0,1,0,0,0,1,1,0,0,0,1,0,0,0) \\
&(1,0,0,0,1,0,0,0,1,0,0,0,1,0,0,0) \\
&(0,1,0,0,0,1,0,0,1,0,0,0,1,0,0,0) \\
&(0,0,1,0,0,0,1,0,1,0,0,0,1,0,0,0) \\
&(0,0,0,0,0,0,0,0,0,1,0,1,1,0,1,0) \\
&(0,0,1,0,0,0,0,1,0,1,0,0,1,0,0,0) \\
&(0,1,0,0,1,0,0,0,0,1,0,0,1,0,0,0)
\end{aligned}
$$
$$
\begin{aligned}
&(1,0,0,0,0,1,0,0,0,1,0,0,1,0,0,0) \\
&(0,0,0,1,0,0,1,0,0,1,0,0,1,0,0,0) \\
&(0,0,0,0,0,0,0,0,1,1,0,0,1,1,0,0) \\
&(0,0,0,0,0,0,0,0,0,0,1,1,1,1,0,0) \\
&(0,1,0,0,0,0,0,1,0,0,1,0,1,0,0,0) \\
&(0,0,1,0,1,0,0,0,0,0,1,0,1,0,0,0) \\
&(0,0,0,1,0,1,0,0,0,0,1,0,1,0,0,0) \\
&(1,0,0,0,0,0,1,0,0,0,1,0,1,0,0,0) \\
&(0,0,0,0,0,0,0,0,1,0,1,0,1,0,1,0) \\
&(0,0,0,0,0,0,0,0,0,1,1,0,1,0,0,1)
\end{aligned}
$$

\subsection{The group $2 \cdot A_4$ with Quaternions and its lattice}
As mentioned previously, Quaternions $\mathbb{H}$ is defined by
$$
\mathbb{H}=\left\{z=x_{0}+x_{1} i+x_{2} j+x_{3} k \mid x_{i} \in \mathbb{R}\right\}
$$
And there is the subalgebra called Hurwitz Quaternionic integers
$$
\mathscr{H}=\{a+i b+j c+k d: a, b, c, d \in \mathbb{Z} \text { or } a, b, c, d \in \mathbb{Z}+1 / 2\} \subset \mathbb{H}
$$
The group generated below whose norm is $1$ is isomorphic to $2 \cdot A_4$.
$$
<i, j, k, \omega>
$$
where $\omega=\frac{1+i+j+k}{2}$. Using the representation defined in the previous section,  the group can also be written as
$$
\left\langle\left(\begin{array}{rrrr}
 & 1 &  &  \\
-1 &  &  &  \\
 &  &  & 1 \\
 &  & -1 &
\end{array}\right),\left(\begin{array}{rrrr}
 &  & 1 &  \\
 &  &  & -1 \\
-1 &  &  &  \\
 & 1 &  &
\end{array}\right),\left(\begin{array}{rrrr}
 &  &  & -1 \\
 &  & -1 &  \\
 & 1 &  &  \\
1 &  &  &
\end{array}\right),
\frac{1}{2}\left(\begin{array}{rrrr}
1 & 1 & 1 & -1 \\
-1 & 1 & -1 & -1 \\
-1 & 1 & 1 & 1 \\
1 & 1 & -1 & 1
\end{array}\right)
\right\rangle
$$
Now consider the inverse map of the bijection defined below for each elements of the group.
$$
\begin{array}{llllll}
\phi: &(\mathbb{Z}\oplus \mathbb{Z}+1/2)^4& \rightarrow& <i,j,k,\omega>&\\
&(c_0,c_1,c_2,c_3)& \mapsto &  c_0 + c_1 i + c_2 j + c_3 k
\end{array}
$$
for example
$$
\begin{array}{lllll}
&\phi^{-1}: &i &\mapsto &(0,1,0,0)\\
&\phi^{-1}: &j &\mapsto &(0,0,1,0)\\
&\phi^{-1}: &k &\mapsto &(0,0,0,1)\\
&\phi^{-1}: &\omega &\mapsto &1/2(1,1,1,1)\\
\end{array}
$$
and the list of the elements is as follows
$$
\begin{array}{lll}
  \pm i\\
  \pm j\\
  \pm k\\
\pm \omega&=\pm \frac{1}{2}(1+i+j+k)&=\mp \omega^4\\
\pm \omega^2&=\pm \frac{1}{2}(-1+i+j+k)&=\mp \omega^5\\
\pm \omega^3&=\pm 1\\
\pm i \omega&=\pm\frac{1}{2}(-1+i-j+k)&=\pm \omega k\\
\pm j \omega&=\pm\frac{1}{2}(-1+i+j-k)&=\pm \omega i\\
\pm k \omega&=\pm\frac{1}{2}(-1-i+j+k)&=\pm \omega j\\
\pm i \omega^{2}&=\pm\frac{1}{2}(-1-i-j+k)&=\pm \omega^{2} j\\
\pm j \omega^{2}&=\pm\frac{1}{2}(-1+i-j-k)&=\pm \omega^{2} k\\
\pm k \omega^{2}&=\pm\frac{1}{2}(-1-i+j-k)&=\pm \omega^{2} i
\end{array}
$$
This yields the following list of $24$ vectors.
$$
\begin{aligned}
&\left(\begin{array}{llll}
1 & 0 & 0 & 0
\end{array}\right)\\
&\left(\begin{array}{llll}
-1 & 0 & 0 & 0
\end{array}\right)\\
&\left(\begin{array}{llll}
0 & 0 & 0 & 1
\end{array}\right)\\
&\left(\begin{array}{llll}
0 & 0 & 0 & -1
\end{array}\right)\\
&\left(\begin{array}{llll}
0 & 0 & 1 & 0
\end{array}\right)\\
&\left(\begin{array}{llll}
0 & 0 & -1 & 0
\end{array}\right)\\
&\left(\begin{array}{llll}
0 & -1 & 0 & 0
\end{array}\right)\\
&\left(\begin{array}{llll}
0 & 1 & 0 & 0
\end{array}\right)
\end{aligned}
$$

$$
\begin{aligned}
&\frac{1}{2}\left(\begin{array}{llll}
-1 & -1 & -1 & -1
\end{array}\right) \\
&\frac{1}{2}\left(\begin{array}{llll}
1 & 1 & 1 & 1
\end{array}\right) \\
&\frac{1}{2}\left(\begin{array}{llll}
1 & -1 & 1 & -1
\end{array}\right) \\
&\frac{1}{2}\left(\begin{array}{llll}
-1 & 1 & -1 & 1
\end{array}\right) \\
&\frac{1}{2}\left(\begin{array}{llll}
1 & 1 & -1 & -1
\end{array}\right) \\
&\frac{1}{2}\left(\begin{array}{llll}
-1 & -1 & 1 & 1
\end{array}\right) \\
&\frac{1}{2}\left(\begin{array}{llll}
-1 & 1 & 1 & -1
\end{array}\right) \\
&\frac{1}{2}\left(\begin{array}{llll}
1 & -1 & -1 & 1
\end{array}\right) \\
&\frac{1}{2}\left(\begin{array}{llll}
-1 & 1 & 1 & 1
\end{array}\right) \\
&\frac{1}{2}\left(\begin{array}{llll}
1 & -1 & -1 & -1
\end{array}\right) \\
&\frac{1}{2}\left(\begin{array}{llll}
-1 & 1 & -1 & -1
\end{array}\right) \\
&\frac{1}{2}\left(\begin{array}{llll}
1 & -1 & 1 & 1
\end{array}\right) \\
&\frac{1}{2}\left(\begin{array}{llll}
-1 & -1 & -1 & 1
\end{array}\right) \\
&\frac{1}{2}\left(\begin{array}{llll}
1 & 1 & 1 & -1
\end{array}\right) \\
&\frac{1}{2}\left(\begin{array}{llll}
1 & 1 & -1 & 1
\end{array}\right) \\
&\frac{1}{2}\left(\begin{array}{llll}
-1 & -1 & 1 & -1
\end{array}\right)
\end{aligned}
$$
And these vectors form the $F_4$ lattice is spanned by the following generator matrix
$$
M(F_4)=\left(\begin{array}{rrrr}
 & 1 & -1 &  \\
 &  & 1 & -1 \\
 &  &  & 1 \\
\frac{1}{2} & -\frac{1}{2} & -\frac{1}{2} & -\frac{1}{2}
\end{array}\right)
$$
can also be spanned by
$$
\left(\begin{array}{cccc}
\frac{1}{2} & \frac{1}{2} & \frac{1}{2} & \frac{1}{2} \\
 & 1 &  &  \\
 &  & 1 &  \\
 &  &  & 1
\end{array}\right)
$$

\subsection{The icosian group order $120$ and $E_8$ Lattice}
I introduced above the subalgebra that is isomorphic to the group $2 \cdot A_4$, actually we can also construct the group $2\cdot A_5$ with quaternions (John Conway, Neil J. A. Sloane[1.P207]) defining $\tau$ and $\sigma$ as
$$
\sigma=\frac{1-\sqrt{5}}{2} \quad \tau=\frac{1+\sqrt{5}}{2}
$$
the $\sigma$ and $\tau$ always satisfy
$$
\begin{array}{l}
\sigma+\tau=-\sigma \tau=1 \\
\sigma^{2}=\sigma+1 \quad \tau^{2}=\tau+1
\end{array}
$$
as they are the answers of the equation $\lambda^2-\lambda-1=0$. Let $x,\omega$ be the following.
$$
x = \frac{1}{2}(i + j\sigma + k\tau) \quad \omega =  \frac{1}{2}(1+i+j+k)
$$
then there is homomorphism such that
$$
\begin{array}{lll}
i & \mapsto &(1,2)(3,4)\\
j &\mapsto &(1,3)(4,2)\\
k & \mapsto &(1,4)(2,3)\\
\omega & \mapsto &(2,3,4)\\
x & \mapsto &(5,2)(3,4)
\end{array}
$$
whose kernel is $\{\pm 1\}$.By simple calculation, you can easily verify the following.
$$
x^4 = \omega^6 = i^4 = j^4 = k^4 = (x\cdot i \cdot \omega)^5=(i\cdot x \cdot \omega^2)^{10}=1
$$
Let $\mathscr{I}$ be the ring formed by icosian, for all the elements $v \in \mathscr{I}$ can be written as
$$
v = \alpha + \beta i + \gamma j + \delta k
$$
where $\alpha,\beta,\gamma,\delta \in \mathbb{Q}[\sqrt{5}]$. Thus it also can be rewritten as
$$
v = (a_1 + a_2 \sqrt{5}) + (a_3 + a_4 \sqrt{5})i + (a_5 + a_6 \sqrt{5})j + (a_7 + a_8 \sqrt{5})k
$$
where $a_i \in \mathbb{Q}$. Now consider the following mapping.
$$
\begin{array}{lllll}
F: &\mathscr{I}&\rightarrow & \mathbb{R}^{4}\\
&v=(a_1 + a_2 \sqrt{5}) + (a_3 + a_4 \sqrt{5})i + \cdots + (a_7 + a_8 \sqrt{5})k&\mapsto &F(v)=(a_1+a_2,a_3+a_4,a_5+a_6,a_7+a_8)
\end{array}
$$
In addition, the same mapping is applied to $\sigma v=(1-\sqrt{5})v/2$, then the lattice formed by the following mapping is isomorphic to $E_8$
$$
\begin{array}{lllll}
&\mathscr{I}^2&\rightarrow & \Lambda_{\mathscr{I}} \subset \mathbb{R}^{8}\\
&(v,\sigma v)&\mapsto &F(v)\oplus F(\sigma v)
\end{array}
$$
whoes the generator matrix is
$$
\frac{1}{2}\left(\begin{array}{llllllll}
1 &  &  & 1 &  & 1 &  & 1 \\
0 & 1 &  & 1 &  &  & 1 & 1 \\
\vdots &  & 1 & 1 &  & 1 & 1 &  \\
 &  &  & 2 &  &  &  &  \\
 &  &  &  & 1 & 1 & 1 & 1 \\
 &  &  &  &  & 2 &  &  \\
\vdots &  &  &  &  &  & 2 &  \\
0 & \cdots &  &  &  & \cdots & 0 & 2
\end{array}\right)
$$
Here is some example, as for $\omega = \frac{1}{2}(1+i+j+k)$
$$
F(\omega) = (\frac{1}{2},\frac{1}{2},\frac{1}{2},\frac{1}{2})
$$
and $\sigma \omega=\frac{1}{4}\{(1-\sqrt{5})+(1-\sqrt{5}) i+(1-\sqrt{5}) j+(1-\sqrt{5}) k\}$, hence
$$
F(\sigma \omega) = (0,0,0,0)
$$
therefore
$$
F(\omega)\oplus F(\sigma \omega)= (\frac{1}{2},\frac{1}{2},\frac{1}{2},\frac{1}{2},0,0,0,0)
$$
For $x$,from the definition
$$
\begin{aligned}
x &=\frac{1}{2}(i+\sigma j+\tau k) \\
&=\frac{1}{2}\left\{i+\left(\frac{1-\sqrt{5}}{2}\right) j+\left(\frac{1+\sqrt{5}}{2}\right) k\right\}
\end{aligned}
$$
Thus
$$
F(x)=(0,\frac{1}{2},0,\frac{1}{2})
$$
And
$$
\sigma x=\left(\frac{1}{4}-\frac{\sqrt{5}}{4}\right) i+\left(\frac{3}{4}-\frac{\sqrt{5}}{4}\right) j-\frac{k}{2}
$$
Thus
$$
F(\sigma x) = (0,0,\frac{1}{2},-\frac{1}{2})
$$
therefore
$$
F(x)\oplus F(\sigma x)=(0,\frac{1}{2},0,\frac{1}{2},0,0,\frac{1}{2},-\frac{1}{2})
$$
In the above way, we obtain the $120$ vectors whose norm is $1$.
\subsection{On the Theta Functions of Lattices}
As one of the simplest example, consider the lattice generated by
$$
M=\left(\begin{array}{ll}
1 & 0 \\
0 & 1
\end{array}\right)=\left(\begin{array}{l}
e_1 \\
e_2
\end{array}\right)
$$
the lattice is simply $\mathbb{Z}^2$
\begin{figure}[h]
 \centering
 \includegraphics[keepaspectratio, scale=0.3]
      {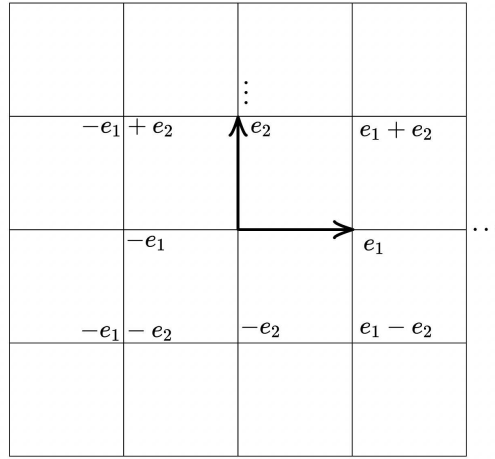}
 \caption{The lattice generated by $M$ which is isomorphic to $\mathbb{Z}^2$}
 \label{}
\end{figure}
The lattice spanned by the generator means that for all elements $x$ in the laatice can be written as the combinations
$$
x = c_1e_1 + c_2e_2 = (c_1 \quad c_2) \cdot M
$$
where $c_i \in \mathbb{Z}$ and let $c=(c_1 \quad c_2)$. As the norm is $N(x)=x \cdot {}^t x=c \cdot M \cdot {}^t (c \cdot M) $, for all norm of the elements of the lattice are written as
$$
N(x) = c \cdot G \cdot {}^tc
$$
where $G$ is gram matrix $G :=M \cdot {}^t M $. In the case of the example, the gram matrix is the identitly matrix. So the norm of all the elements in the lattice can be written as
$$
N(x)= c\cdot {}^t c = c^2_1 + c^2_2
$$
functions with the same number of norms of the lattice as coefficients can be defined as follows
$$
\Theta(q) = \sum_{x \in \Lambda} q^{N(x)}
$$
where $\Lambda$ is the lattice spanned by the generators, and we call the function is theta function of lattices.In this case, the function will be
$$
\Theta(q) = \sum_{x \in \Lambda} q^{N(x)} = \sum^{\infty} _{c_1 =-\infty} \sum^{\infty} _{c_2 =-\infty} q^{c^2_1 + c^2_2}
$$
using Jacobi theta functions
$$
\begin{aligned}
\theta_{2}(q)&=\sum_{n=-\infty}^{\infty} q^{(m+1 / 2)^{2}} \\
\theta_{3}(q)&=\sum_{n=-\infty}^{\infty} q^{m^{2}} \\
\theta_{4}(q)&=\sum_{n=-\infty}^{\infty}(-q)^{m^{2}}
\end{aligned}
$$
the theta function of the lattice is rewritten as
$$
\Theta(q) = \theta_3(q)^2=1 + 4q + 4q^2 + 4q^4 + 8q^5 +  4q^8  + 4q^9 + 8q^{10} +\cdots
$$
from the definition, the coefficients are equal to the numbers of coordinates on the same circle in the following figure.
\begin{figure}[h]
 \centering
 \includegraphics[keepaspectratio, scale=0.15]
      {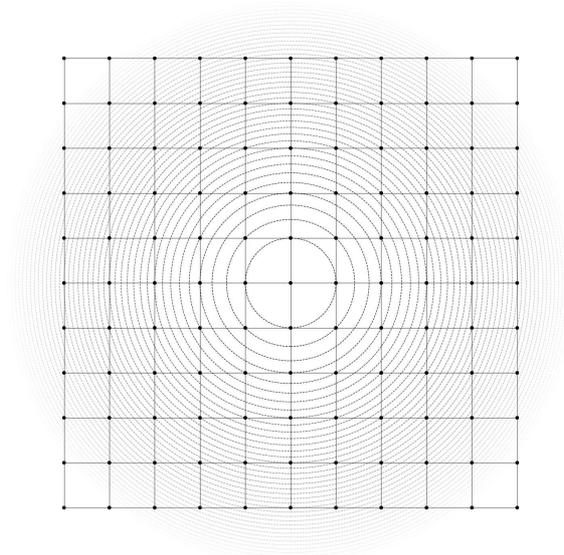}
 \caption{Relation between vectors of the same distance.}
 \label{}
\end{figure}

\end{document}